\documentclass{amsart}
\begin{document}

\vfuzz2pt 
\hfuzz2pt 
\newtheorem{thm}{Theorem}[section]
\newtheorem{corollary}[thm]{Corollary}
\newtheorem{lemma}[thm]{Lemma}
\newtheorem{proposition}[thm]{Proposition}
\newtheorem{defn}[thm]{Definition}
\newtheorem{remark}[thm]{Remark}
\newtheorem{example}[thm]{Example}
\newtheorem{fact}[thm]{Fact}
\
\newcommand{\norm}[1]{\left\Vert#1\right\Vert}
\newcommand{\abs}[1]{\left\vert#1\right\vert}
\newcommand{\set}[1]{\left\{#1\right\}}
\newcommand{\Real}{\mathbb R}
\newcommand{\eps}{\varepsilon}
\newcommand{\To}{\longrightarrow}
\newcommand{\BX}{\mathbf{B}(X)}
\newcommand{\A}{\mathcal{A}}
\newcommand{\onabla}{\overline{\nabla}}
\newcommand{\hnabla}{\hat{\nabla}}
\newcommand{\f}{\mathbf{f}}

\def\proof{\medskip Proof.\ }
\font\lasek=lasy10 \chardef\kwadrat="32 
\def\kwadracik{{\lasek\kwadrat}}
\def\koniec{\hfill\lower 2pt\hbox{\kwadracik}\medskip}

\newcommand*{\C}{\mathbf{C}}
\newcommand*{\R}{\mathbf{R}}
\newcommand*{\Z}{\mathbf {Z}}

\def\sb{f:M\longrightarrow \C ^n}
\def\det{\hbox{\rm det}\, }
\def\detc{\hbox{\rm det }_{\C}}
\def\i{\hbox{\rm i}}
\def\tr{\hbox{\rm tr}\, }
\def\rk{\hbox{\rm rk}\,}
\def\vol{\hbox{\rm vol}\,}
\def\Im {\hbox{\rm Im}\, }
\def\Re{\hbox{\rm Re}\, }
\def\interior{\hbox{\rm int}\, }
\def\e{\hbox{\rm e}}
\def\pu{\partial _u}
\def\pv{\partial _v}
\def\pui{\partial _{u_i}}
\def\puj{\partial _{u_j}}
\def\puk{\partial {u_k}}
\def\div{\hbox{\rm div}\,}
\def\Ric{\hbox{\rm Ric}\,}
\def\r#1{(\ref{#1})}
\def\ker{\hbox{\rm ker}\,}
\def\im{\hbox{\rm im}\, }
\def\I{\hbox{\rm I}\,}
\def\id{\hbox{\rm id}\,}
\def\exp{\hbox{{\rm exp}^{\tilde\nabla}}\.}
\def\cka{{\mathcal C}^{k,a}}
\def\ckplusja{{\mathcal C}^{k+1,a}}
\def\cja{{\mathcal C}^{1,a}}
\def\cda{{\mathcal C}^{2,a}}
\def\cta{{\mathcal C}^{3,a}}
\def\c0a{{\mathcal C}^{0,a}}
\def\f0{{\mathcal F}^{0}}
\def\fnj{{\mathcal F}^{n-1}}
\def\fn{{\mathcal F}^{n}}
\def\fnd{{\mathcal F}^{n-2}}
\def\Hn{{\mathcal H}^n}
\def\Hnj{{\mathcal H}^{n-1}}
\def\emb{\mathcal C^{\infty}_{emb}(M,N)}
\def\M{\mathcal M}
\def\Ef{\mathcal E _f}
\def\Eg{\mathcal E _g}
\def\Nf{\mathcal N _f}
\def\Ng{\mathcal N _g}
\def\Tf{\mathcal T _f}
\def\Tg{\mathcal T _g}
\def\diff{{\mathcal Diff}^{\infty}(M)}
\def\embM{\mathcal C^{\infty}_{emb}(M,M)}
\def\U1f{{\mathcal U}^1 _f}
\def\Uf{{\mathcal U} _f}
\def\Ug{{\mathcal U} _g}
\def\[f]{{\mathcal U}^1 _{[f]}}
\def\hnu{\hat\nu}
\def\gnu{\nu_g}
\title{Bochner's technique for statistical structures}
\author{Barbara Opozda}

\subjclass{ Primary: 53B05, 53C05, 53A15, 53B20}

\keywords{affine connection, statistical structure, curvature
tensors, Laplacian, Hedge's theory,  Bochner's technique}

\thanks{The research supported by the NCN grant UMO-2013/11/B/ST1/02889
and a grant of the TU  Berlin} \maketitle

\address{Instytut Matematyki UJ, ul. \L ojasiewicza  6, 30-348 Cracow,
Poland}

\email{Barbara.Opozda@im.uj.edu.pl}

\vskip 1in \noindent

\vskip 0.5in \noindent \begin{abstract}  The main aim of this paper
is to extend Bochner's technique to statistical structures. Other
topics related to this technique are also introduced to the theory
of statistical structures. It deals, in particular,  with Hodge's
theory,  Bochner-Weitzenb\"ock and Simon's type formulas. Moreover,
a few global and local theorems on the geometry of statistical
structures are proved, for instance, theorems saying that under some
topological and geometrical conditions a statistical structure must
be trivial. We also introduce a new concept of sectional curvature
depending on statistical connections. On the base of this notion we
study the curvature operator and prove some analogues of well-known
theorems from Riemannian geometry.
\end{abstract}

\maketitle

\vskip1truecm

\section{Introduction}
 The main tool of  the Bochner technique is the Levi-Civita
connection. Our  purpose is to show that the technique can be
extended to the class of statistical  connections.

We shall study the following four cases:
\newline
i) A torsion-free connection $\nabla$ is statistical for a metric
tensor field $g$, that is, $\nabla g$ is symmetric.
 A statistical structure, that is, a pair  $(g,\nabla)$, where $\nabla$ is  statistical for $g$
  is also called a Codazzi pair.
\newline
ii)  A statistical  structure $(g,\nabla)$ is equiaffine, that is,
there is a volume form $\nu$ such that $\nabla\nu =0$.
\newline
iii) A statistical structure $(g,\nabla)$  is equiaffine relative to
the volume form $\gnu$ determined by $g$. It is equivalent to the
condition $\tr _g\nabla g(X, \cdot,\cdot)=0$ for every $X$. We shall
call such structures trace-free.
\newline
iv) For a statistical structure $(g,\nabla)$ the curvature tensors
for $\nabla$ and its conjugate connection $\overline \nabla$ are the
same.

The oldest examples of statistical structures are the induced
structures (consisting of the second fundamental form and the
induced connection) on locally strongly convex hypersurfaces in
$\R^{n+1}$
 endowed with an equiaffine transversal vector field (in other words -- with relative normalization).
Within the theory of equiaffine hypersurfaces the case iii)
corresponds to  Blaschke hypersurfaces,
  the  case iv) -- to equiaffine  spheres.

 But the
majority of statistical structures is outside the class of
hypersurfaces. Even using the structures obtained on hypersurfaces
one can easily modify them and get structures which are not
realizable on hypersurfaces. For instance, the product of equiaffine
ovaloids is equipped with the product statistical structure but it
cannot be realized as a
 locally strongly convex hypersurface in  any $\R^{n+1}$.
 It is also easy to find examples of statistical structures which
  are non-realizable on hypersurfaces even locally.

In Section 2 we  provide preliminary information on divergences for
volume forms and connections and establish few  integral formulas
useful in proving classical Bochner's theorems and their
generalizations to statistical structures (e.g. in Section 10).

Basic notions for statistical structures and  their subclasses
listed above are introduced and discussed in Section 3. Examples of
various types of statistical structures are given in Section 4.

In Section 5 we prove a few global theorems saying that under some
topological and geometrical assumptions two statistical structures
with the same metric must be identical or a statistical structure
must be trivial, that is, the statistical connection is the
Levi-Civita connection.

On a statistical manifold one can define  various Laplacians (acting
on differential forms). First we have the Laplacian for the
underlying Riemannian structure. One can also define the
codifferential $\delta ^\nabla$ relative to a statistical connection
$\nabla$ and then set
$$\Delta ^\nabla =\delta ^{\onabla}d +d\delta ^{\onabla},$$
where $\onabla$ is the conjugate connection for $\nabla$. For this
Laplacian we prove basic properties for compact manifolds and
 Hodge-type theorems (Sections 8, 9).
 A differential form $\omega$ will
  be called $\nabla$-harmonic if $\Delta ^\nabla\omega=0$.
  Bochner's technique for vector fields and harmonic 1-forms
is developed in Sections 10 and 11.  Bochner-Weitzenb\"ock formulas
for Laplacians acting on differential forms are computed in Section
11. There we also compute Simons'type formulas for the Laplacians of
the square of the length of any tensor field. The
Bochner-Weitzenb\"ock curvature operator can be also applied to
other tensor fields, in particular, to the metric tensor field of
statistical structures.

Another aim  of the paper is to introduce  a notion of  sectional
curvature for statistical structures. The curvature tensor of
$\nabla$ does not have, in general,
 as good symmetries as the curvature tensor of the Levi-Civita connection.
 But one can modify it and get a tensor field  with  the same symmetries as the Riemannian curvature tensor.
  Using the modified tensor one can define  an appropriate notion of  sectional curvature
  and the corresponding curvature operator.
   After the modification there is still a problem with the second Bianchi identity,
   which plays an essential role in  many  theorems, e.g. Schur's lemma or Tachibana's theorem.
   In the general case of statistical structures, Schur's lemma does not hold.
  By restricting  considerations to the class iv) we get appropriate
analogues of  such theorems.

\section{Divergences and integral formulas}
All the objects considered in this paper  are of class $\mathcal C
^\infty$. All connections are linear and torsion-free.  Let $M$ be
an $n$-dimensional manifold with a fixed volume form $\nu$. For any
vector field $X$ on $M$ its Lie derivative $\mathcal L_X\nu$ is an
$n$-form, hence
\begin{equation}
\mathcal L_X\nu=(\div ^\nu X)\nu.
\end{equation}
 The function $\div ^\nu X$ is the divergence
relative to the volume form $\nu$. A  divergence can also be defined
relative to a  connection. Namely, if $\nabla$ is a  connection,
then
\begin{equation}\div ^\nabla X=\tr \{Y\to \nabla_YX\}\end{equation}
for a vector field $X$. More generally, for any
tensor field $s$ of type $(1,k)$ we have
\begin{equation}
(\div ^\nabla s )(X_1,...,X_k)=\tr \{Y\to (\nabla
_Ys)(X_1,...,X_k)\}.
\end{equation}


\begin{lemma} For any  connection $\nabla$ on $M$  and a tensor field  $S$ of type
$(1,1)$ we have
\begin{equation}\label{lemat_ze_sladem}
X\tr S=tr\nabla _XS.
\end{equation}
\end{lemma}
\proof     In the equality (\ref{lemat_ze_sladem}) both sides depend
on $X$ in a tensorial way. Let $x\in M$ and $X\in T_xM$. Take a
local frame $e_1,...,e_n$ and its dual frame $\theta _1,...,
\theta_n$ defined around $x$ and such that $\nabla e_i=0$, $\nabla
\theta^i=0$  at $x$ for  $i=1,...n$. Since at $x$
$$0=(\nabla _X\theta ^i)(Se_i)=X(\theta ^i(Se_i))-\theta ^i(\nabla _X(Se_i)),$$
we have (at $x$)
\begin{eqnarray*}
X\tr S&=&\sum_{i=1}^nX(\theta ^i(Se_i)) =\sum _{i=1}^n \theta ^i(\nabla _X(Se_i))\\
&=& \sum_{i=1}^n\theta ^i((\nabla _XS)e_i)=\tr \nabla _X S.
\end{eqnarray*}\koniec

The equality (\ref{lemat_ze_sladem}) can be written as
\begin{equation}
\nabla_X(\tr S)=\tr (\nabla _XS).
\end{equation}
Denote by  $R$  the curvature tensor of  $\nabla$ and by
$Ric$ its Ricci tensor.
\begin{lemma}\label{second_lemma}
Let $\nabla$ be a  connection on $M$. For a vector field $X$ on $M$
we set
\begin{equation}
S_XY=\nabla_YX.
\end{equation}
$S_X$ is a $(1,1)$-tensor field and $\div ^{\nabla}X=\tr \, S_X$. We
have
\begin{equation}\label{about_L_S}
\mathcal L_X =\nabla_X -S_X.
\end{equation}
 For any vector fields  $X, Y$
on $M$ the following formula holds
\begin{equation}\label{lematRic(X,Y)}
\begin{array}{lcr}
\div ^\nabla(\nabla _XY)&=& Ric (X,Y) +\tr(\nabla _XS_Y) +\tr (S_Y\circ S_X)\\
&=& Ric(X,Y) +\ \ X(\tr S_Y) +\tr (S_Y\circ S_X)\\
&=& Ric(X,Y) +X(\div^\nabla Y) +\tr (S_Y\circ S_X).
\end{array}
\end{equation}
\end{lemma}
\proof  The equality (\ref{about_L_S}) is well known and it immediately follows from the
fact that $\mathcal L_X$ and $\nabla _X$ are differentiations. For
any vector fields $X,Y, Z$ on $M$ the following equalities hold
\begin{eqnarray*}
\nabla _Z\nabla _XY &=& R(Z,X)Y +\nabla _X\nabla _ZY +\nabla _{[Z,X]}Y\\
&=& R(Z,X)Y +\left( \nabla _X \nabla _Z Y-\nabla_{\nabla _XZ}Y\right) +\nabla _{\nabla _ZX}Y\\
&=& R(Z,X)Y +(\nabla _XS_Y)Z +(S_Y\circ S_X)Z.
\end{eqnarray*}
Taking the trace relative to $Z$ on both sides we obtain the
required equality.\koniec

\begin{lemma}\label{divergences}
Let $\nabla$ be a connection and $\nu$ be a volume form on $M$. Then
\begin{equation}
\div ^\nabla X =  \div ^\nu X+\tau (X)
\end{equation}
for any $X\in \mathcal X(M)$, where $\nabla _X\nu =-\tau (X)\nu$.
\end{lemma}
\proof It follows from the equality $ \nabla_X\nu=\mathcal L _X\nu
+S_X\nu.$\koniec

According to \cite{NS}, by an equiaffine structure we mean a pair
$(\nabla ,\nu)$ consisting of a connection  $\nabla$  and a volume
form $\nu$  such that
 $\nabla \nu=0$. Thus for an equiaffine structure $(\nabla, \nu)$ we
 have
 $\div ^\nabla = \div ^\nu$.

\begin{lemma}\label{th_int_Ric(X,Y)}
If  $M$ is a compact manifold with a  volume form $\nu$ and $\nabla$
is a  connection on $M$,
 then for any vector fields
  $X,Y$ on $M$  we have
\begin{equation}\label{prebasic}
\int _MRic(X,Y)\nu =\int_M \tr S_Y (\div ^{\nu} X)\nu -\int_M \tr
(S_Y\circ S_X)\nu+ \int_M\tau(\nabla _XY)\nu.
\end{equation}
\newline
In particular, if $\nabla \nu =0$ then
\begin{equation}\label{basic}
 \int _MRic(X,Y)\nu =\int_M \tr S_Y \tr S_Y \, \nu-\int_M \tr (S_Y\circ S_X)\nu.
\end{equation}
\end{lemma}
\proof Let  $X,Y$ be  vector fields on  $M$. Set $\varphi =\tr S_Y$.
We  have
$$ \mathcal L_X(\varphi\nu )=(X\varphi )\nu +\varphi \mathcal L_X\nu= (X(\tr S_Y))\nu+
(\tr S_Y \div ^\nu X)\nu.
$$ and by Stokes' theorem we get
$$ \int _M X(\tr S_Y)\nu =-\int _M \tr S_Y (\div ^\nu X)\nu .$$
Using this equality, Lemmas \ref{second_lemma}, \ref{divergences}
and the divergence theorem we obtain the result.\koniec

In the 2-dimensional case we get

\begin{proposition}
Let  $M$ be a compact  $2$-dimensional manifold with an equiaffine
structure $(\nabla ,\nu)$.
 For each vector field  $X$ on  $M$
 we have
 \begin{equation}
\int _MRic(X,X)\nu=2\int _M (\det S_X )\nu .
 \end{equation}
\end{proposition}
\proof  For any endomorphism  $A$ of a 2-dimensional  vector space
we have $(\tr A)^2 -\tr A^2 =2\det A$\koniec

If an endomorphism $A$
of a real vector space is diagonalizable, then $\tr (A^2) \ge
0$.
 Therefore, by (\ref{basic}), we get
\begin{proposition}\label{calka_Ric(X,X)}
If  $M$  is a compact manifold with an equiaffine structure $(\nabla
,\nu)$,  a vector field  $X\in\mathcal X (M)$
 is without divergence and  $S_X$ is diagonalizable at each point  $M$,
 then
\begin{equation}
\int _MRic(X,X) \nu \le 0.
\end{equation}
\end{proposition}
\medskip

\section{Statistical structures}
For a tensor field $s$  and a connection $\nabla$ the notation
$\nabla s (X,...)$ will stand for $(\nabla _Xs)(...)$.

Let $g$ be a positive definite Riemannian tensor field on a manifold $M$. We assume that $M$ is oriented.
 Denote by $\hnabla$ the
Levi-Civita connection for $g$  and by $\gnu$ the volume form
determined by $g$. 
 We shall study (torsion-free) connections
$\nabla$ satisfying the following Codazzi condition:

\begin{equation}\label{symmetry}
(\nabla_X g)(Y,Z)=(\nabla _Yg)(X,Z)
\end{equation}
for all $X,Y,Z\in T_x M$, $x\in M$. A structure  $(g,\nabla)$
satisfying (\ref{symmetry}) is called a statistical structure. We
shall call a connection $\nabla$ satisfying (\ref{symmetry}) a
statistical connection for $g$.  Since
\begin{equation}\label{nabla_Xnu_g}2\nabla_X\nu _g=\tr _g (\nabla _Xg)(\cdot,\cdot
)\nu _g,\end{equation} the condition $ \nabla \gnu =0$
 is equivalent to the condition
\begin{equation}\label{trace-free}
\tr _g(\nabla _Xg)(\cdot,\cdot )=0
\end{equation}
for every $X\in TM$. If $\nabla$ is statistical for $g$ and  (\ref{trace-free}) is satisfied, we
shall say that the statistical structure is trace-free.

For any connection $\nabla$ one defines its conjugate $\onabla$
relative to $g$ by the formula
\begin{equation}
g(\nabla _XY,Z)+g(Y,\onabla _XZ)=Xg(Y,Z).
\end{equation}
It is known that if $(g,\nabla)$ is trace-free then so is
$(g,\onabla)$, if  $(g,\nabla)$ is a statistical structure then so
is $(g,\onabla)$, see e.g. \cite{NS}. Hence statistical structures
go in pairs.

If $R$ is the curvature
tensor for $\nabla$ and $\overline R$ is the curvature tensor for
$\onabla$, then we have, \cite{NS},
\begin{equation}\label{R_and_oR}
g(R(X,Y)Z,W)=-g(\overline{R}(X,Y)W,Z).
\end{equation}
It follows that
\begin{equation}\label{ricci_for_both}
\overline{Ric} (Y,W)=-\tr _gg(R(\cdot,Y)\cdot, W),
\end{equation}
where $\overline{Ric}$ is the Ricci tensor of $\onabla$. The
function
\begin{equation}
\rho=\tr _gRic(\cdot,\cdot)
\end{equation}
will be called the scalar curvature of $(g,\nabla)$. Similarly we
define the scalar curvature $\overline\rho$ for $(g,\onabla)$ and we
have the usual scalar curvature $\hat\rho$ for $(g,\hnabla)$. By
(\ref{ricci_for_both}) we have
\begin{equation}\label{scalar_curvatures}
\rho =\overline\rho .
\end{equation}
From now on in this section we assume that $\nabla$ is statistical
for $g$. If $K$ is the difference tensor between $\nabla$ and
$\hnabla$, that is,
\begin{equation}
\nabla _XY=\hnabla _XY+K_XY,
\end{equation}
then
\begin{equation}
\onabla_XY=\hnabla _XY-K_XY
\end{equation}
and
\begin{equation}
\hnabla _XY={1\over 2}(\nabla_XY +\onabla_XY).
\end{equation}
 $K(X,Y)$ will stand for $K_XY$. Since $\nabla$  and $\hnabla$ are  without
torsion, $K$  as a $(1,2)$-tensor is symmetric. We have $ (\nabla
_Xg)(Y,Z)=(K _Xg)(Y,Z)=-g(K_XY,Z)-g(Y,K_XZ)$. It is now clear that
the symmetry of $\nabla g$ and $K$ implies the symmetry of $K_X$
relative to $g$ for each $X$. The converse also holds. Namely, if
$K_X$ is symmetric relative to $g$ then we have
\begin{equation}\label{nablag_gK}
\nabla g(X,Y,Z) = -2g(K_XY,Z).
\end{equation}
Set
\begin{equation}
 E=\tr _gK(\cdot,\cdot).
\end{equation}
If $\tau (X):=\tr K_X$  then $\tau (X)=g(E,X)$. By (\ref{nablag_gK})
we have
\begin{equation}\label{tr_g_nablag_tau}
\tr _g \nabla g(\cdot, \cdot, Z) =-2g(E,Z)=-2\tau (Z).\end{equation}
Comparing this equality with (\ref{nabla_Xnu_g}) we see that $\nabla
_Z\nu _g=-\tau (Z)\nu _g$ (compare also with  Lemma
\ref{divergences}). We have $g(\nabla _Xg,\nabla_X g)=4g(K_X,K_X)$
and since $\onabla g(X,Y,Z)=2g(K_XY,Z),$ we also have
\begin{equation}\label{g(K,K)}
g(\onabla_Xg,\onabla_Xg)=4g(K_X,K_X)=g(\nabla _Xg,\nabla _Xg).
\end{equation}
Consequently
\begin{equation}
g(\onabla g,\onabla g)=4g(K,K) =g(\nabla g,\nabla g).
\end{equation}
Observe also that
\begin{equation}
g(\nabla X,\onabla X) =g(\hnabla X,\hnabla X) -g(K_X,K_X)
\end{equation}
for any $X\in \mathcal X(M)$. Indeed, one has
\begin{eqnarray*}g(\nabla _YX,\onabla _YX)&=&g(\hnabla _YX+K_YX,\hnabla
_YX-K_YX)\\&&=g(\hnabla_YX,\hnabla _YX)-g(K_YX,K_YX).\end{eqnarray*}
Similarly
\begin{equation}\label{g(nabla,nabla)+g(onabla,onabla)}
g(\nabla X,\nabla X) +g(\onabla X,\onabla X) =2\{ g(\hnabla
X,\hnabla X) +g(K_X,K_X)\}.
\end{equation}
For a vector field $X$ we have three  $(1,1)$-tensor fields $S_X$,
$\hat S_X$ and $\overline S_X$ defined  by $S_XY=\nabla _YX$, $\hat
S_XY=\hnabla _YX$, $\overline S_XY=\onabla _YX$. It is clear that
\begin{equation}\label{trSXoSX}
\tr S_X=\tr \hat S_X +\tau (X),\ \ \ \ \ \ \tr \overline S_X=\tr \hat S_X-\tau(X)
\end{equation}
for every $X\in \mathcal X(M)$. It is known that
\begin{equation}\label{from_Nomizu_Sasaki}
R(X,Y)=\hat R(X,Y) +(\hnabla_XK)_Y-(\hnabla_YK)_X+[K_X,K_Y].
\end{equation}
Writing the same equality for $\onabla$ and adding both equalities we
get
\begin{equation}\label{R+oR}
R(X,Y)+\overline R(X,Y) =2\hat R(X,Y) +2[K_X,K_Y].
\end{equation}
We also have
 \begin{eqnarray*}
\tr\{X \to [K_X,K_Y]Z\}&=& \sum _{i}^n g(e_i,[K_{e_i},K_Y]Z)\\
&&\ \ \ \ = \sum_i^n
(g(e_i,K_{e_i}K_YZ))-g(e_i,K_YK_{e_i}Z))\\
&& \ \ \ \ \ \ \ \ \ \ \ \ \ \ =\sum_i^n (g(K_{e_i} e_i,K_YZ) -g(K_Ye_i, K_Ze_i))\\
&&\ \ \ \ \ \ \ \ \ \ \ \ \ \ \ \ \ \ \ \ \ \ \ \ \ \ \ \ \ =-g(K_Y,K_Z) +\tau (K(Y,Z)).
\end{eqnarray*}
Choose now a point $x_0$ and an orthonormal frame $e_1,...,e_n$
around $x_0$ such that $\hnabla e_i=0$ for $i=1,...n$ at $x_0$.
Having $Y,Z\in T_{x_0}M$ we extend the vectors to vector fields, say
$Y,Z$, around $x_0$ in such a way that $\hnabla Y=\hnabla Z=0$ at
$x_0$. In particular, $[Y,Z]=0$  at $x_0$.
 We obtain at $x_0$

\begin{eqnarray*}
\sum_{i=1}^n&&[ g((\hnabla_{e_i}K) (Y,Z), e_i) - g( (\hnabla _YK)
(e_i, Z), e_i)]\\
&& =(div ^{\hnabla} K)(Y,Z)
-\sum_{i=1}^nYg(K(e_i,Z), e_i)\\
&&=(div ^{\hnabla} K)(Y,Z)
-\sum_{i=1}^nYg(K(e_i,e_i), Z)\\
&&=(div ^{\hnabla} K)(Y,Z)- Yg(E,Z)\\
&&=(div ^{\hnabla} K)(Y,Z)- Y\tau(Z)\\
&&=(div ^{\hnabla} K)(Y,Z)- \hnabla \tau (Y,Z).
\end{eqnarray*}
Therefore
\begin{equation}\label{Ricci_tensor}
Ric(Y,Z)=\widehat{Ric}(Y,Z) + (div^{\hnabla}K)(Y,Z) -\hnabla
\tau(Y,Z) +\tau (K(Y,Z)) -g(K_Y,K_Z).
\end{equation}
 It follows that
\begin{equation}\label{Ric+oRic}
Ric (Y,Z)+\overline{Ric} (Y,Z)=2\widehat{Ric}(Y,Z) -2g(K_Y,K_Z) +2\tau(K(Y,Z)).
\end{equation}
In particular, if $(g, \nabla )$ is trace-free then
\begin{equation}2\widehat{Ric}(X,X)\ge Ric(X,X)+\overline{Ric}(X,X).\end{equation}
The above formulas also yield
\begin{equation}\label{symetria_Ric}
Ric(Y,Z)-Ric (Z,Y)=-g((\hnabla K(Y,e_i,Z),e_i)+ g((\hnabla
K(Z,e_i,Y),e_i)= - d\tau (Y,Z).
\end{equation}
Hence $\nabla$ is Ricci-symmetric if and only if $d\tau=0$. The
following lemma follows from  formulas (\ref{R_and_oR}),
(\ref{from_Nomizu_Sasaki}) and (\ref{R+oR}).

\begin{lemma}\label{przeniesiony_lemat}
Let $(g,\nabla) $ be a statistical structure. The following
conditions are equivalent:
\newline
{\rm 1)} $R=\overline R$,
\newline
{\rm 2)} $\hnabla K$ is symmetric,
\newline
{\rm 3)} $g(R(X,Y)Z,W)$ is skew-symmetric relative to $Z,W$.
\end{lemma}
From 2) and (\ref{symetria_Ric}) we see that the condition
$R=\overline R$ implies the symmetry of $Ric$.
 We have proved
\begin{proposition}Let $(g,\nabla)$ be a statistical structure.
$Ric$ is symmetric if and only if $d\tau=0$. If $R=\overline R$ then
$Ric=\overline{Ric}$ is symmetric.
\end{proposition}


Taking now the trace  relative to $g$ on both sides of (\ref{Ric+oRic}) and taking into
account that $\rho =\overline\rho$, we get
\begin{equation}\label{theorema_egregium}
\hat\rho =\rho +|K|^2- |E|^2.
\end{equation}
In the case where $\nabla$ is the induced connection on a Blaschke
hypersurface in $\R ^{n+1}$ and $g$ is the Blaschke metric, the
equality (\ref{theorema_egregium}) (with $E=0$) is known as the
affine  theorema egregium. Indeed, if $H$  is the affine mean
curvature then  $H=n^2(n-1)\rho$ and $|K|^2=4n(n-1)J$, where $J$ is
the Pick invariant.

For an orthonormal frame $e_1,...,e_n$ we have
$$ |K|^2=g(K,K)=\sum  _{i,j,k}g(K_{e_i}e_j,e_k)^2,\ \ \ \ \ \  |E|^2=g(E,E)=\sum _{j,k} g(K_{e_j}e_j,e_k)^2.
$$
Thus  $|K|^2-|E|^2\ge 0$ on $M$. If $|K|=|E|$  then
$0=g(K_{e_i}{e_j} ,e_k) =g(K_{e_k}e_i,e_j)$ for every $k$ and $i\ne
j$. It follows that $K_X$ is a multiple of the identity for each
$X$, which is possible only for $K=0$. Thus we have

\begin{proposition}
The functional
$$\mathfrak{scal}:\{ statistical\ connections\ for \ g\} \ni\nabla \to \tr _g Ric\in \mathcal C^{\infty}(M)$$
attains its maximum for the Levi-Civita connection at each point of
$M$. Conversely, if $\nabla$ is a statistical connection  for $g$
and $\mathfrak{scal}$ attains its maximum for $\nabla$ at each point
on $M$, then $\nabla$ is the Levi-Civita connection for $g$.
\end{proposition}

\begin{corollary}
Let $(g,\nabla)$ be a statistical structure on $M$ and $\rho\ge
\hat\rho$ on $M$. Then $\nabla$ is the Levi-Civita connection for
$g$.
\end{corollary}

 We shall also study equiaffine statistical structures. By an
 equiaffine statistical structure on $M$ we mean a triple $(g,\nabla
 ,\nu)$, where $(g,\nabla)$ is a statistical structure and $\nu$ is
 a volume form on $M$ such that $\nabla \nu =0$. Let us emphasize
 that $\nu $ is not necessarily the volume form $\nu _g$.

\medskip

\medskip
\section{Examples } The theory of affine
hypersurfaces in $\R ^{n+1}$ is a natural  source of statistical
structures. For the theory we refer to \cite{LSZ} or\cite{NS}. We
recall here only some basic facts.

 Let $\mathbf{f} :M\to \R
^{n+1}$ be a locally strongly convex  hypersurface.  For simplicity
assume that $M$ is connected and oriented. Let $\xi$ be a
transversal vector field on $M$. We define the induced  volume form
$\nu _\xi$ on $M$ (compatible with the given orientation) as follows
$$\nu_\xi (X_1,...,X_n)=\det (\mathbf {f}_*X_1,...,\mathbf{f}_*X_n,\xi).$$
We  also have the induced connection $\nabla$ and the second
fundamental form $g$ defined by the Gauss formula:
$$D_X\mathbf{f}_*Y=\mathbf{f}_*\nabla _XY +g(X,Y)\xi,$$
where $D$ is the standard flat connection on $\R ^{n+1}$. Since the
hypersurface is locally strongly convex, $g$ is  definite. By
multiplying $\xi$ by $-1$, if necessary, we can assume that $g$ is
positive definite. A transversal vector field is called equiaffine
if $\nabla \nu_\xi=0$. This condition is equivalent to the fact that
$\nabla g$ is symmetric, i.e. $(g,\nabla)$ is a statistical
structure. It means, in particular, that for a statistical structure
obtained on a hypersurface by a choice of a transversal vector
field, the Ricci tensor of $\nabla$ is automatically symmetric.

For later use  recall the notion of the shape operator and the Gauss
equations. Having a chosen equiaffine transversal vector field  and
differentiating it we get the Weingarten formula
$$D_X\xi= -\mathbf{f}_*\mathcal SX.$$
The tensor field $\mathcal S$ is called the shape operator for
$\xi$. If $R$ is the curvature tensor for the induced connection
$\nabla$ then
\begin{equation}\label{Gauss_equation_for R}
R(X,Y)Z=g(Y,Z)\mathcal SX-g(X,Z)\mathcal SY.
\end{equation}
This is the Gauss equation for $R$. The Gauss equation for
$\overline R$ is the following
\begin{equation}\label{Gauss_equation_for_oR}
\overline R(X,Y)Z=g(Y,\mathcal SZ)Y-g(X,\mathcal SZ)X.
\end{equation}
In particular, the dual connection is projectively flat. Recall also
that the form $g(\mathcal SX,Y)$ is symmetric for any equiaffine
transversal vector field.

 For a locally strongly convex hypersurface
there are  infinitely many equiaffine transversal vector fields. In
fact, if $\xi$ is any equiaffine transversal vector field for
$\mathbf{f}$ (for instance a metric normal vector field) and $\phi$
is a nowhere vanishing function on $M$, then
$\tilde\xi=\mathbf{f}_*Z+\phi\xi$ is equiaffine, where $g(Z,X)=
X\phi$.  We  also have  the volume form determined by $g$ on $M$. In
general, this volume form is not covariant constant relative to
$\nabla$. It can be proved that there is a unique equiaffine
transversal vector field $\xi$ such that $\nu_\xi =\nu_g$. This
unique transversal vector field is called the affine normal vector
field. The second fundamental form for the affine normal is called
the Blaschke metric. If the affine lines determined by the affine
normal vector field meet at one point or are parallel then the
hypersurface is called an affine sphere. In the first case the
sphere is called proper in the second one  improper. The class of
affine spheres is very large. There exist a lot of conditions
characterizing affine spheres. For instance, a hypersurface is an
affine sphere if and only if $R=\overline R$, see Lemma
\ref{sfery_afiniczne} below.

As we have already observed, if $\nabla$ is a connection on a
hypersurface induced by an equiaffine transversal vector field then
the conjugate connection $\onabla$ is projectively flat. Therefore
the projective flatness of the conjugate connection is a necessary
condition for $(g,\nabla)$ to be realizable as the induced structure
on a hypersurface.

We will now make few remarks on  statistical structures in general, that is,
 possibly non-realizable on hypersurfaces.

As we have mentioned in the introduction, the cartesian product of
statistical manifolds is a statistical manifold which cannot be
realized as a locally strictly convex hypersurface.

We shall now produce other statistical structures which are
non-realizable on hypersurfaces.

 Observe that if $(g,\nabla)$ is a statistical
structure with the difference tensor $K$ and $\phi$ is any smooth
function on $M$ then
$$\tilde\nabla=\nabla +\phi K=\hnabla +(1+\phi)K$$
is a statistical connection for $g$. Moreover, if $(g,\nabla)$ is
trace-free then so is $(g,\tilde \nabla)$. If $\phi$ is constant
and $R=\overline R$ then $\overline{\tilde R}=\tilde R$. Indeed, in
this case we  have $\hnabla ((1+\phi)K)=(1+\phi)(\hnabla K)$ and we
can  now use the above lemma. We now have
\begin{proposition} Assume that $\mathbf{f}:M\to \R ^{n+1}$, where $n>2$, be a locally
strongly convex affine sphere  equipped with the statistical
structure $(g,\nabla)$ described above. Assume that the sectional
curvature for $g$ is not constant on $M$. There is no $t\in
\R\setminus \{0,-2\}$ such that $(g,\tilde\nabla)$ is realizable on
a hypersurface $ \tilde {\mathbf{ f}}:M\to\R^{n+1}$, where $\tilde\nabla =\nabla
+tK$.
\end{proposition}
\proof The connection $\onabla$ is projectively flat, hence
$$R(X,Y)Z=\overline {R}(X,Y)Z=\gamma (Y,Z)X-\gamma (X,Z)Y$$ for some
$(0,2)$-tensor field (the normalized Ricci tensor for $\onabla$).
Since $\hnabla K$ is symmetric, we have
$$R(X,Y)=\hat R(X,Y) +[K_X,K_Y].$$ We now
have
$$\tilde R(X,Y) =\hat R(X,Y) +(1+t)^2[K_X,K_Y]= R(X,Y)
+t(2+t)[K_X,K_Y].$$ Suppose that $\overline{\tilde R}$ ($=\tilde R$)
is projectively flat.  If $t\ne 0$ and $t\ne -2$ then
$$[K_X,K_Y]Z=\gamma _1(Y,Z)X -\gamma _1(X,Z)Y$$
for some $(0,2)$-tensor field $\gamma_1$. But it means that
$\hnabla$ is projectively flat, which contradicts the assumption
that the sectional curvature of $g$ is not constant.\koniec

Note that all affine spheres whose  Blaschke metric has constant
sectional curvature are known. These are quadrics (for which $\nabla
=\hnabla$) or hypersurfaces given by the equations
$$x_1\cdot\cdot\cdot x_{n+1}=c,$$
where $x_1,....,x_{n+1}$ are the canonical coordinates in $\R
^{n+1}$ and $c=const\ne 0$, see Theorem 2.2.3.18 in \cite{LSZ}.

Come back to the observation that for a statistical structure
realizable on a hypersurface in $\R^{n+1}$ the Ricci tensor of its
connection must be symmetric. Assume we have  a locally strongly
convex hypersurface equipped with an equiaffine transversal vector
fiels and the induced statistical structure on it. Assume that it is
not trace-free. We have the non-zero vector field $E$ and its dual
form $\tau$. Let $\phi$ be a function on $M$ such that $d\phi\ne
\tau$. Consider the connection $\tilde \nabla=\nabla +\phi K$. Then
$\tilde \tau =(1+\phi)\tau $. Since the Ricci tensor for the
statistical structure is symmetric if and only if $d\tau\equiv 0$
and $d(\phi\tau)=d\phi\wedge \tau$, we see that $d\tilde \tau\ne 0$.
Hence $\tilde\nabla$ is not Ricci symmetric and consequently
$(g,\tilde\nabla)$ cannot be realized on a hypersurface.

\bigskip

\section{Further properties of statistical structures}
For a given metric tensor $g$ one has, in general, many statistical
connections. Given a connection one  also has, in general, many
metric tensor fields constituting with the connection a statistical
structure. But if we impose additional conditions on the structures
and manifolds, the situation might change drastically.

Consider, for instance,  the following problem. Let $(g_1,\nabla)$
be a trace-free statistical structure on $M$. Does there exist
another metric tensor $g_2$ (non-homothetic to $g_1$)  on $M$ for
which $(g_2, \nabla )$ is a trace-free statistical structure. For
structures realizable on hypersurfaces at least $3$-dimensional the
answer is negative. The answer is also negative for $2$-dimensional
ovaloids in $\R ^3$. In the last case, the theorem is, in fact, true
for abstract compact $2$-dimensional manifolds of genus $0$. To
illustrate this type of consideration we give a proof of this fact.

Assume  that $M$ is $2$-dimensional connected and oriented. Having a
statistical structure on $M$ and, in particular, a positive definite
metric tensor field, we also have the underlying complex structure
on $M$. Denote by $\mathfrak{S}_k(M)$, $k>1$, the space of all
symmetric $k$-covariant  tensor fields $s$ for which
$$\tr _gs (\cdot,\cdot, X_3,..., X_k)=0.$$
Let $z=x+\i y$ be an isothermal coordinate on $M$ and
$X=\frac{\partial}{\partial x}$, $Y=\frac{\partial}{\partial y}$.
Denote by $\mathfrak{ S}_k^{\mathbf C}(M)$ the space of all complex
symmetric $k$-forms on $M$. By a straightforward purely algebraic
computation one gets
\begin{lemma} For $s\in \mathfrak{ S} _k(M)$ the symmetric complex form
$$\Phi(s)= [s(X,...,X)-\i s(Y,X,...,X)]dz^k$$
is well-defined on the whole of $M$, i.e. it is independent of a
choice of isothermal coordinates. The mapping
$$\Phi: \mathfrak{ S}_k(M)) \ni s\to
 \Phi (s)\in \mathfrak{ S}_k^{\mathbf C}(M)$$
 is a linear isomorphism
(over $\mathcal C^\infty (M)$).
\end{lemma}

A symmetric  tensor $s$ is called a Codazzi tensor for a connection
$\nabla$ if $\nabla s$ is symmetric. By computing the Cauchy-Riemann
or, in a general version, Vekua-Carleman equations one gets

\begin{lemma}\label{carleman}
If  a symmetric tensor field $s\in \mathfrak S_k(M)$ is a Codazzi
tensor  for the Levi-Civita connection $\hnabla$ then the form
$\Phi(s)= [s(X,...,X)-\i s(Y,X,...,X)]dz^k$ is holomorphic. If $s$
is Codazzi for any torsion-free connection then the form $\Phi (s)$
is pseudo-holomorphic.
\end{lemma}
By the Riemann-Roch theorem (or the index method) one  knows that if
$M$ is  compact then a pseudo-holomorphic symmetric complex $k$-form
is either constantly zero or its  zeros  are isolated and their
number (counted with multiplicities) is equal to $-2\chi (M)$, where
$\chi(M)$ is the Euler characteristic of $M$.

We   can now prove  the following rigidity result due  to U. Simon,
\cite{S}.

\begin{thm}
Let $M$ be a connected compact surface of genus $0$. If
$(g_1,\nabla)$, $(g_2,\nabla)$ are two  trace-free structures on $M$
then $g_1=cg_2$ on $M$ for some constant number $c$.
\end{thm}
\proof We can assume that $M$ is oriented. Since $\nabla \nu
_{g_1}=\nabla \nu _{g_2}=0$, by multiplying $g_2$ by a constant we
can assume that $\nu _{g_1}=\nu _{g_2}$.

Define $g=g_1+g_2$ and $h=g_1-g_2$. Since both statistical
structures are trace-free, we have that the cubic forms $\nabla g$ i
$\nabla h$ are symmetric. Observe that $\tr _gh(\cdot, \cdot)=0$.
There is a basis $e_1,e_2$ of $T_pM$ which is $g_1$-- orthonormal
and such that $g_2(e_i,e_j)=\lambda_i \delta_{ij}$. By the
assumption $\nu _{g_1}=\nu _{g_2}$ we have $\lambda _1\lambda _2=1$.
The vectors $\frac{e_1}{\sqrt{1+\lambda_1}},
\frac{e_2}{\sqrt{1+\lambda_2}}$ form a $g$--orthonormal basis. We
now have
$$\tr_gh=\frac{1-\lambda _1}{1+\lambda _1} +\frac{1-\lambda_2}{1+\lambda
_2}=0.$$

Using Lemma \ref{carleman} for $h$  and then the Riemann-Roch
theorem  finishes the proof.\koniec
\bigskip

The same consideration as in the above proof can be applied to
surfaces of other topological types. For instance, we have

\begin{thm}
Let $M$ be a connected compact surface of genus $1$. If
$(g_1,\nabla)$, $(g_2,\nabla)$ are two  trace-free structures on $M$
and $g_1=cg_2$
 at one point of $M$ then $g_1=cg_2$ on the whole $M$.
\end{thm}

 Other results typical for affine differential geometry are those
saying when the induced connection must be the Levi-Civita
connection for the second fundamental form. A similar problem is
interesting in the case of abstract statistical structures. For
instance, using the above considerations we obtain

\begin{thm}
Let $M$ be a  connected compact surface and $(g, \nabla)$ be a
trace-free statistical structure on $M$. Let $R=\overline R$. If $M$
is of genus $0$ then $\nabla=\hnabla$ on $M$. If $M$ is of genus $1$
and $K=0$ at one point of $M$ then $\nabla =\hnabla$ on $M$.
\end{thm}
\proof It is sufficient to consider the symmetric cubic form
$C(X_1,X_2,X_3)= g(K(X_1, X_2), X_3)$. Since $R=\overline R$ implies
that $\hnabla C$ is symmetric, we have that the complex form
$$[C(X,X,X)-\i C(Y,X,X)]dz^2$$
is holomorphic. Using the Riemann-Roch theorem finishes the
proof.\koniec

In the higher-dimensional case we have the following theorem.
 If  a metric
tensor field $g$ is given, $\hat k(X\wedge Y)$ will denote the
sectional curvature by the plane spanned by $X,Y$ if these vectors
are linearly independent.
\begin{thm}
Let $M$ be a compact manifold equipped with a trace-free statistical
structure $(g,\nabla)$ such that $R=\overline R$. If the sectional
curvature $\hat k$ for $g$ is positive then $\nabla=\hat \nabla$.
\end{thm}

\proof Since $R=\overline R$, $\hnabla K$ is symmetric. Consider the
function on the unit sphere bundle $UM$
$$\alpha: UM\ni V\to g(K(V,V),V)\in \R.$$
Let $e_1$ be a point on $UM$, where $\alpha$ attains its maximum.
 Denote the maximal value by $\lambda
_1$. Let $u\in U_pM$ be orthogonal to $e_1$, that is, $u$ is tangent
to $U_pM$ at $e_1$. Take the curve $\beta (t)=\cos te_1 +\sin t u$.
Since $\alpha$ attains a maximum at $e_1$, by differentiating
$\alpha\circ\beta$ at $t=0$ we obtain
\begin{equation}
g(K(e_1,e_1),u)=0,\ \ \ \ 2g(K(e_1,u),u)-g(K(e_1,e_1),e_1)\le 0.
\end{equation}
The first formula yields $K(e_1,e_1)=\lambda _1e_1$ for some
$\lambda _1$, that is, $e_1$ is an eigenvector for $K_{e_1}$. Let
$e_1,...e_n$ be an orthonormal eigenbasis for $K_{e_1}$ and  let
$\lambda _2$,..., $\lambda _n$ be eigenvalues corresponding to
$e_2,..., e_n$ respectively. By the above inequality we have
\begin{equation}
\lambda _1- 2\lambda _i\ge 0
\end{equation}
for $i=2,...,n$.

 For any $u\in U_pM$, take the $\hnabla$-geodesic  $\gamma$ in $M$ with $\gamma
(0)=p$ and $\gamma '(0)=u$. Let $e_1(t)$ be  the vector field along
$\gamma$ obtained by the  parallel displacement relative to
$\hnabla$ of the vector $e_1$. We get a vector field $e_1(t)$. Since
$\alpha$ attains the maximum at $e_1$, by differentiating
$\alpha\circ\gamma$ we obtain
$$0= \frac{d}{dt}_{|t=0}g(K(e_1(t),e_1(t)),e_1(t))=g(\hnabla K
(u,e_1,e_1),e_1)$$ and
\begin{equation}\label{nierownosc}
0\ge \frac{d^2}{dt^2}_{|t=0}g(K(e_1(t),e_1(t)),e_1(t))=g(\hnabla
^2K(u,u,e_1,e_1), e_1).\end{equation}
 We have
\begin{eqnarray*}
&&( \hnabla ^2 K)(X,Y,Z,W)-(\hnabla ^2 K)(Y,X,Z,W)= (\hat
R(X,Y)\cdot
K)(Z,W)\\
&&\ \ \ \ \  \ =\hat R (X,Y)(K(Z,W)) -K(\hat R(X,Y)Z,W)-K(Z,\hat
R(X,Y)W)
\end{eqnarray*}
for every vectors $X,Y,Z,W$. Therefore
\begin{eqnarray*}
&&g( (\hnabla ^2 K)(U,V,U,V),V)-g((\hnabla ^2 K)(V,U,U,V),V)\\
&&\ \ \ \ \  \ =g(\hat R (U,V)(K(U,V)),V) -g(K(\hat
R(U,V)U,V),V)-g(K(U,\hat R(U,V)V),V)
\end{eqnarray*}
for any $U,V$. Using also the symmetry of $\hnabla K$ we get
\begin{eqnarray*}
&&g((\hnabla ^2K)(e_i,e_i, e_1,e_1),e_1)=g((\hnabla
^2K)(e_1,e_1,e_i,e_i),
e_1)\\
&&\ \ \ \ +g(\hat R(e_i,e_1)(K(e_i, e_1)),e_1) -g(K(\hat R
(e_i,e_1)e_i,e_1), e_1)-g(K(e_i,\hat R(e_i,e_1)e_1),e_1)
\end{eqnarray*}
for every $i=1,...,n$. Using now  (\ref{nierownosc}), the symmetries
of $K$  and the assumption that $\tr _gK=0$ we see that
\begin{eqnarray*}
&&0\ge \sum _{i=1}^n g((\hnabla ^2K)(e_i,e_i, e_1,e_1), e_1)\\
&&\ \ \  =\sum _{i=2}^n -\lambda_ig(e_i, \hat R (e_i,e_1)e_1)
+\lambda_1\hat k(e_i\wedge e_1)-\lambda_i \hat k(e_i\wedge e_1)\\
&&\ \ \ = \sum_{i=2}^n \hat k (e_1\wedge e_i)(\lambda_1-2\lambda
_i).
\end{eqnarray*}
Since the curvature $\hat k >0$ and $\lambda _1-2\lambda _i\ge 0$
for $i=2,...,n$, we have $\lambda _1-2\lambda _i=0$ for $i=2,..,n$.
Using now  the assumption that $\tr K_{e_1}=0$ we see that
$\lambda_1=0$. Hence $\alpha=0$ on $M$ and consequently $K\equiv 0$
on $M$.\koniec


If $(M,g)$ is a compact oriented  Riemannian manifold and $UM$
denotes the unit sphere bundle then for every tensor field $s$ of
type $(0,k)$ we have
\begin{equation}\label{RosI}
\int _{UM} (\hnabla s) (U,..., U)dU=0
\end{equation}
and
\begin{equation}\label{RosII}
\int _{UM}\tr _g(\hnabla s)(\cdot,\cdot, U,...,U)dU=0,
\end{equation}
where
$$\int _{UM}fdU=\int _{x\in M}\left(\int
_{U_xM}f\nu_{g_x}\right)\nu_g$$ for any continuous  function
$f:UM\to \R$.
 These Ros' integral  formulas can be found in \cite{R} and \cite{MRU}.

The formula (\ref{RosII}) adapted to $(1,k)$-tensor fields says the
following. If $s$ is a tensor field of type $(1,k)$ then
\begin{equation}\label{RosII'}
\int_{UM} (div^{\hnabla}s)(U,...,U)dU=0.
\end{equation}
We can now prove

\begin{thm} Let $M$ be a compact oriented manifold and
$(g,\nabla)$ be a statistical structure on it. Then
\begin{equation}
\int_{UM}Ric (U,U)dU=\int_{UM}\overline {Ric} (U,U)dU
\end{equation}
and
\begin{equation}\label{w_twierdzeniu_o_dwoch_liczbach}
\begin{array}{rcl}
&&\int_{UM}Ric (U,U)dU =\int_{UM}\widehat{Ric}(U,U)dU-\int
_{UM}g(K_U,K_U)dU\\
&&\ \ \ \ \ \ \ \ \ \ \ \ \ \ \ \  \ \ \ \ \ \ \ \ \ \ \ \ \ \ \ \ \
\ \ \ \ \ \ \ \ \ \ \ \ \ \ \ \ \ \ \ \ \ \  \ \ \ \ \ \ \ \
 +\int _{UM}\tau(K(U,U))dU.\end{array}
\end{equation}
In particular, if $(g,\nabla)$ is trace-free then the following
equality of two numbers
$$\int _{UM}Ric
(U,U)dU=\int_{UM}\widehat{Ric}(U,U)dU$$ implies that $\nabla
=\hnabla$ on $M$.
\end{thm}
\proof The first two formulas follow from (\ref{Ricci_tensor}),
(\ref{Ric+oRic}) and Ros' integral formulas applied to $K$ and
$\tau$. To prove the last assertion it is now sufficient to use
(\ref{w_twierdzeniu_o_dwoch_liczbach}). \koniec


\begin{corollary}
Let  $M$ be an ovaloid in $\R^{n+1}$ equipped with an equiaffine
transversal vector field and $g$ be the corresponding second
fundamental form, $\nabla$   the induced connection and $\mathcal S$
the shape operator. Then
\begin{equation}
n\int_{UM}g(\mathcal SU,U)dU=vol (S^{n-1})\int_M \tr \mathcal S ,
\end{equation}
where $vol (S^{n-1})$ is the volume of the unit sphere in the
standard Euclidean space $\R ^n$.
\end{corollary}
\proof By formulas (\ref{Gauss_equation_for R}) and
(\ref{Gauss_equation_for_oR}) one has
$$Ric (Y,Z)=g(Y,Z)\tr \mathcal S -g(Z,\mathcal SY),$$
$$\overline {Ric}(Y,Z)=(n-1)g(\mathcal SZ,Y).$$
Using the above theorem we get the assertion.\koniec

In particular, if the third fundamental form $g(\mathcal
S\cdot,\cdot)$ is positive definite (like in the case where the
transversal vector field is a metric normal for a locally strictly
convex hypersurface) then $\int _M\tr \mathcal S>0$.

It is known that if $g$ is the second fundamental form on a
hypersurface in $\R ^{n+1}$  corresponding to a transversal vector
field and $\nabla ^2g=0$ then $\nabla =\hnabla$. In the compact
case one has a stronger result.

\begin{proposition}
Let $M$ be  compact oriented and $(g,\nabla)$ be a statistical
structure on $M$. Then
\begin{equation}
\int_{UM}\nabla ^2g(U,U,U,U)dU=6\int_{UM}\Vert K(U,U)\Vert ^2dU
\end{equation}
In particular, if $\int _{UM}\nabla ^2g(U,U,U,U)dU=0$ then $\nabla
=\hnabla$.
\end{proposition}
\proof  Denote by $C$ the symmetric cubic form $\nabla g$. We have
$$\nabla ^2 g(U,U,U,U)=\hnabla C (U,U,U,U) +(K_UC)(U,U,U)$$
and
\begin{eqnarray*}
&&(K_UC)(U,U,U)=-3C(K_UU,U,U)\\&&\ \ \ \ \ =-3\sum_{i=1}^n
g(K(U,U),e_i)C(e_i,U,U)=6\sum _{i=1}^ng(K(U,U),e_i)^2.
\end{eqnarray*}
Using the integral formula (\ref{RosI}) finishes the proof.\koniec

\bigskip

\section{Laplacians for statistical structures}
We adopt the following convention
\begin{eqnarray*}
d\omega
(X_0,...,X_k)&= &\sum _{i=0}^k (-1)^iX_i(\omega(X_0,...,\hat
X_i,...,X_k))\\
&& -\sum_{0\le i<j\le k}(-1)^{i+j}\omega([X_i,X_j],X_0,...,\hat
X_i,...,\hat X_j,...,X_k)
\end{eqnarray*}
and consequently
\begin{equation}\label{d_omega_connection}
d\omega (X_0,...,X_k)=\sum _{i=0}^k (-1)^i(\nabla
_{X_i}\omega)(X_0,...,\hat X_i,...,X_k)
\end{equation}
for any torsion-free connection $\nabla$.

Let $M$ be an oriented manifold and $g$  a positive definite metric tensor field on $M$.
 We have the standard Hodge Laplacian

\begin{equation}
\Delta =\delta d+d\delta.
\end{equation}
If  $f\in \mathcal C^{\infty}(M)$ then
\begin{equation}
\Delta f=-\div^{\gnu} grad\, f.
\end{equation}
 If $(g,\nabla) $ is a trace-free statistical structure and $f$ is a function
then
$$\Delta f=-\tr \{Y\longrightarrow\nabla _Y(grad\, f)\}.$$
 For any connection
$\nabla$ we can define the Hessian of a function $f$
$$Hess^{\nabla}f(X,Y) =\nabla ^2 f(X,Y)=\nabla_X(\nabla_Yf)=X(df(Y))-df(\nabla_XY).$$
The Hessian is a tensor field of type $(0,2)$. If the connection
$\nabla$ is torsion-free, the Hessian $Hess^\nabla f$  is symmetric.
 We have
\begin{lemma}\label{BW_for_functions}
For any statistical structure $(g,\nabla)$ and any  function $f\in \mathcal C ^{\infty}(M)$  we have
$$\Delta f=-\tr _g Hess^{\nabla}f(\cdot,\cdot)-df(E).$$
\end{lemma}
\proof We have
\begin{eqnarray*}
Hess^{\nabla} \, f(X,Y)& =& X(df( Y )) -df(\nabla _XY) \\
&=& X(df(Y)) -df(\hnabla _XY)-df(K(X,Y))\\
 &=&  Hess ^{\hnabla}f(X,Y) - df(K(X,Y)).
\end{eqnarray*}

\koniec

From now on we assume that $(g,\nabla)$ is a statistical structure
and we shall use the notions introduced in Section 3. For a
statistical structure $(g,\nabla)$ we shall study  a Laplacian
relative to the connection $\nabla$. Note that the Laplacian
 we  propose is different than the
 Laplacian relative to a connection (called also the Lichnerowicz Laplacian) defined   as
  $\tr _g (\nabla ^2s)(\cdot, \cdot )$ for a tensor field $s$.

  If
$f$ is a function then one sets
\begin{equation}\label{Delta_nabla_f}
\Delta^\nabla f=-\div^\nabla grad\, f.\end{equation} It is clear
that
$$\Delta^\nabla f =-\tr _g (\onabla _{\cdot }df)(\cdot)=-\tr _g
(\onabla ^2f).$$
 This Laplacian acting on functions was introduced and studied in \cite{W}. We shall
extend it to the operator acting on  differential forms. We define
the codifferential relative to $\nabla$  acting on differential
forms as follows
\begin{equation}
\delta ^\nabla \omega =-\tr _g\nabla \omega(\cdot,\cdot,...).
\end{equation}
A differential form $\omega$ will be called $\nabla$-coclosed if
$\delta^\nabla\omega=0$. We shall use
\begin{lemma}\label{lemma_about_K}
For any $k$-form $\omega$ and any orthonormal frame $e_1,...,e_n$ we
have
  $$\sum_{i} (K_{e_i}\omega)(e_i,...)=-\iota_E\omega,$$
  where $\iota$  stands for  the interior product.
\end{lemma}
\proof We have
$$\sum _{i=1}^n(K_{e_i}\omega)(e_i,Y_2,...,Y_k)=-\sum_{i=1}^n \omega
(K_{e_i}e_i,Y_2,...,Y_k)-\sum _{i=1}^n\sum_{l=2}^k\omega
(e_i,Y_2,...,K_{e_i}Y_l,...,Y_k).$$ If $\omega$ is a $1$-form the
second term on the right hand side  does not appear. If $\omega$ is
of degree at least 2, we fix $2\le l\le k$. Using the symmetries of
$K$ we compute
\begin{eqnarray*}
\sum _{i=1}^n\omega
(e_i,Y_2,...,K_{e_i}Y_l,...,Y_k)&=&\sum_{i,j=1}^n\omega
(e_i,Y_2,...,g(K_{e_i}Y_l,e_j)e_j,...,Y_k)\\
&=& \sum_{i,j=1}^ng(K_{e_i}Y_l,e_j)\omega(e_i,Y_2,...,e_j,...,Y_k)\\
&=& \sum_{i<j}g(K_{e_i}Y_l,e_j)\omega(e_i,Y_2,...,e_j,...,Y_k)\\
&& +\sum _{j<i}g(K_{e_i}Y_l,e_j)\omega(e_i,Y_2,...,e_j,...,Y_k)\\
&=&  \sum_{i<j}g(K_{e_i}Y_l,e_j)\omega(e_i,Y_2,...,e_j,...,Y_k)\\
&& - \sum_{i>j}g(K_{e_j}Y_l,e_i)\omega(e_j,Y_2,...,e_i,...,Y_k)\\
&=&0.
\end{eqnarray*}
\koniec

From the above proof we also have
\begin{lemma}\label{in_the_long_proof}
Let $\alpha$ be an $(m+1)$-form, $m\ge 1$, and an index $l$, $1\le
l\le m$,  be fixed. We have
\begin{equation}
\sum_{j=1}^n\alpha (e_j, X_1,..., K _{e_j}X_l,..., X_m)=0.
\end{equation}
\end{lemma}
 Extending the definition (\ref{Delta_nabla_f}) we set
$$\Delta  ^\nabla \omega =\delta ^{\onabla} d\omega
+d\delta^{\onabla}\omega$$ for any differential form $\omega$.
 Immediate consequences of Lemma \ref{lemma_about_K} and the classical Weitzenb\"ck formula:
  $\delta\omega=-\tr _g\hnabla \omega (\cdot,\cdot,...)$ are the
following relations
\begin{equation}\label{relation_codifferentials}
\delta=\delta ^\nabla -\iota _E=\delta^{\onabla} +\iota _E,
\end{equation}
\begin{equation}\label{relation_laplacians}
\Delta ^\nabla =\Delta -\mathcal L_E.
\end{equation}
\begin{lemma}\label{delta_nabla_delta_nabla}
For any statistical connection $\nabla$ we
 have $\delta^\nabla\delta ^\nabla =0$.
\end{lemma}
\proof By (\ref{relation_codifferentials}) it is sufficient to
observe that $\iota_E \delta  +\delta \iota _E=0$. The equality
trivially holds  for $0$-  and $1$-forms. Let $\omega$ be a
$k$-form, where $k\ge 2$. Take an orthonormal frame $e_1,..,e_n$
around a fixed point $x_0\in M$ such that $\hnabla e_i=0$ at $x_0$.
Extend vectors $X_1,...,X_{k-2}\in T_{x_0}M$  to local vector fields
around $x_0$ in such a way that $\hnabla X_j=0$ at $x_0$ for
$j=1,...,k-2$.
Using the standard Weitzenb\"ock formula we get at $x_0$
\begin{eqnarray*}
(\iota _E\delta \omega)(X_1,...,X_{k-2})&=&-\sum_{i=1}^n e_i(\omega (e_i,E,X_1,...,X_{k-2})),\\
(\delta \iota_E\omega)(X_1,...,X_{k-2})&=& -\sum_{i=1}^n
e_i(\omega(E,e_i, X_1,...,X_{k-2})).
\end{eqnarray*}\koniec

Denote by $\mathcal F ^k(M)$ the space of all (smooth) differential
forms of degree $k$ and by $\mathcal F(M)$ the algebra of all
differential forms on $M$. The metric tensor $g$  extended to the
bundle of tensors on $TM$ will be denoted by $g$.

By Lemma \ref{delta_nabla_delta_nabla} the operator $\delta^\nabla$
determines  the exact sequence
$$0\leftarrow\mathcal F^0(M)\leftarrow\mathcal F^1(M)\leftarrow
...\ .$$

We shall say that a differential form $\omega$ is $\nabla$-harmonic
if $\Delta^\nabla\omega =0$.


Assume now that $(g,\nabla,\nu)$ is an equiaffine statistical
structure. Let $\nu=\varphi \gnu$. We have
$$\hnabla _X\nu =d\varphi (X)\gnu , \ \ \ \hnabla _X\nu =-K_X \nu
=-K_X (\varphi \gnu ) =\varphi(\tr K_X)\gnu.$$ Thus $$\tau =d\, log
\, \varphi ,$$ where $\tau (X)=\tr K_X$. For any orthonormal frame
$e_1,...,e_n$ and any form $\omega$ we now have
\begin{eqnarray*}
\delta (\varphi\omega)&=& -\tr _g\hnabla
(\varphi\omega)(\cdot,\cdot,...)=\varphi\delta\omega -\sum_{i=1}
^n\varphi\tau(e_i)\omega(e_i,...)\\
&=&\varphi \delta\omega - \varphi \omega(\sum _{i=1} ^n
\tau(e_i)e_i,...)= \varphi\delta\omega -\varphi\iota _E\omega
.\end{eqnarray*} Thus
\begin{equation}\label{codifferential_phi_omega}
\delta(\varphi\omega)=\varphi\delta^{\onabla}\omega .
\end{equation}

\begin{lemma}\label{d_delta_onabla_adjoint}
Let $M$ be a compact manifold equipped with an equiaffine
statistical structure $(g,\nabla,\nu)$. For any differential forms
$\omega$, $\eta$ we have
\begin{equation} \int _Mg(\omega,d\eta)\nu
=\int _M g(\delta ^{\onabla} \omega,\eta)\nu.\end{equation}
\end{lemma}
\proof Using (\ref{codifferential_phi_omega}) we obtain
\begin{eqnarray*}
\int_M g(\omega,d\eta)\nu &=&\int _M
g(\omega,d\eta)\varphi\gnu=\int_M
g(\varphi\omega,d\eta)\gnu\\
&=& \int _M  g(\delta (\varphi \omega),\eta)\gnu=\int _M \varphi
g(\delta ^{\onabla }\omega,\eta)\gnu=\int_M g(\delta ^{\onabla},
\eta )\nu .
\end{eqnarray*} \koniec

\begin{corollary}\label{Delta_nabla_self_conjugate}
 For an equiaffine statistical structure $(g,\nabla,\nu)$ on a compact manifold
$M$
\begin{equation}
\int _Mg(\Delta ^\nabla \omega ,\eta)\nu=\int _M g(\omega, \Delta
^\nabla \eta)\nu=\int_M g(d\omega, d\eta)\nu +\int _M
g(\delta^{\onabla}\omega,\delta ^{\onabla}\eta)\nu .
\end{equation}

\end{corollary}

\begin{corollary}Let $M$ be a compact manifold equipped with an equiaffine
statistical structure $(g,\nabla, \nu)$. A differential form is
$\nabla$-harmonic if and only if it is closed and
$\onabla$-coclosed.
\end{corollary}
If a differential form $\omega$ is $\onabla$-parallel then $\omega$
is $\onabla$-coclosed. It is also closed, because $\onabla$ is
torsion-free. Therefore we have
\begin{proposition}
For any statistical structure $(g,\nabla)$ $\onabla$-parallel forms
are $\nabla$-harmonic.
\end{proposition}

\bigskip

\section{Formal adjoint  operators for statistical connections on compact manifolds}

Assume that $(g,\nabla)$ is a statistical structure. We shall
construct an appropriate formal adjoint operators for $\nabla$ using
a standard procedure.

We shall  use the musical isomorphism notation, that is,
 $X^{\flat}(Y)=g(X,Y)$,
$g(\alpha ^{\sharp},Y)=\alpha (Y)$ for $\alpha \in T_xM^*,\ X,Y\in
T_xM,\ x\in M.$

 We restrict our consideration to  sections of tensor bundles.
 Let $\mathcal E$ be a vector subbundle of the tensor bundle over
$M$. The set of all smooth sections of this bundle will be denoted
by $\mathcal C^{\infty}(M\leftarrow \mathcal E)$. The metric tensor
and the connections $\nabla$, $\overline\nabla$ are extended to the
bundle $\mathcal E$. If $s_1, s_2$ are sections of $\mathcal E$ then
$$Xg(s_1,s_2)=g(\nabla _Xs_1,s_2) +g(s_1,\onabla _X s_2)$$
for any $X\in TM$. We now take the bundle $$HOM (TM,\mathcal E)=\sum
_{x\in M}HOM (T_xM,\mathcal E_x),$$ where $HOM (T_xM,\mathcal E_x)$
is the space of all linear mappings from $T_xM$ to $\mathcal E_x$.

A section of the bundle $HOM(TM,\mathcal E)$ can be treated as a
mapping which sends a vector field $X\in \mathcal X(M)$ to a section
of $\mathcal E$. The metric tensor $g$ is again extended to this
bundle and as usual denoted by the same letter $g$. The same deals
with the extended connections $\nabla$ and $\onabla$. If $S$ is a
section of $HOM (TM,\mathcal E)$ then for any $X,Y\in \mathcal X(M)$
we have

$$(\nabla _XS)Y= \nabla _X(SY)-S(\nabla _XY).$$
We now regard  $\nabla$ as  a mapping
$$\nabla : \mathcal C ^{\infty}(M\leftarrow \mathcal E)\longrightarrow
\mathcal C ^{\infty}(M\leftarrow HOM(TM,\mathcal E)),$$  where
$$\nabla s= \{ \mathcal X(M)\ni X \to \nabla _X s\in \mathcal C
^{\infty}(M\leftarrow \mathcal E)\},$$ that is,
$$(\nabla s)(X)=\nabla _X s.$$
 If $\tilde \nabla$ is a connection on $M$ (possibly different
 than $\nabla$) then
\begin{equation}\label{nabla^2}
\begin{array}{rcl}
(\tilde\nabla _Y(\nabla s))(X) &= &\tilde\nabla
_Y((\nabla s)(X))-(\nabla s) (\tilde\nabla_YX)\\
&=& \tilde\nabla _Y(\nabla _Xs)-\nabla _{\tilde\nabla _YX}s.
\end{array}
\end{equation}
In particular, if $\tilde\nabla=\nabla$ then
 $\nabla^2=\nabla (\nabla s)$.
We denote $(\nabla _Y(\nabla s))(X)$ by $\nabla ^2 _{Y,X} s$. If $s$
is a section of  $\mathcal E$  and $S$ is a section of $HOM(TM,
\mathcal E)$ then $g(\nabla s, S)$ is a function on $M$. In order
to compute it take a local orthonormal frame  $e_1,..., e_n$ in a
neighborhood of a fixed point  $x_0\in M$  such that $\hnabla e_j=0$
at  $x_0$. In particular, $\sum _{i=1}^n \nabla _{e_i}e_i=E$ and
$\sum _{i=1}^n \onabla _{e_i}e_i=-E$ at $x_0$. We obtain at $x_0$
\begin{equation}\label{iloczyn_skalarny_nabla_sS}
\begin{array}{rcl}
 g(\nabla s, S) &=& \sum
_{i=1}^ng( \nabla _{e_i}s, S{e_i})\\
&=&\sum _{i=1}^n {e_i}g( s,Se_i) -\sum _{i=1}^n g( s, \onabla _{e_i}
(Se_i))\\
&=&( \div ^{\gnu} \alpha ^\sharp) -\sum _{i=1}^n g( s, \onabla
S(e_i,e_i))+g(s,SE),
\end{array}
\end{equation}
where $\alpha$ is a  $1$-form on  $M$ given by $\alpha (X)=g( s, SX)
$. Indeed, we have the following equalities at $x_0$
$$\div^{\gnu} \alpha ^{\sharp}=\sum _{i=1}^n g( \hnabla _{e_i}\alpha
^{\sharp}, e_i)=\sum _{i=1}^n( e_i g(\alpha ^{\sharp},e_i) -g(
\alpha ^\sharp, \hnabla _{e_i}e_i))=\sum_{i=1}^ne_i(\alpha (e_i))
$$
and
$$\sum _{i=1}^n\onabla _{e_i}(Se_i)=\sum_{i=1}^n\onabla S(e_i,e_i) +\sum _{i=1}^nS(\onabla
_{e_i}e_i)=\sum_{i=1}^n\onabla S(e_i,e_i)-SE.$$ Similarly we have at
$x_0$
\begin{eqnarray*}
\div ^\nabla\alpha ^\sharp =\sum_{i=1}^n(e_ig(\alpha^\sharp
,e_i)-g(\alpha ^\sharp ,\onabla _{e_i}e_i))&=&\sum _{i=1}^n
e_i(\alpha (e_i))+\alpha(E)\\&=& \sum _{i=1}^n
e_i(g(s,Se_i))+g(s,SE)\end{eqnarray*}
 and consequently
\begin{equation}\label{iloczyn_skalarny_nabla_sS_nabla}
g(\nabla s,S)=\div ^\nabla \alpha ^\sharp-\sum_{i=1}^n g(s,\onabla
S(e_i,e_i)).
\end{equation}

Assume  that $M$ is compact and oriented. We first consider the
scalar products (both denoted by $\langle ,\rangle$) on the infinite
dimensional vector spaces $\mathcal C ^{\infty}(M\leftarrow \mathcal
E)$, $\mathcal C ^{\infty}(M\leftarrow HOM(TM,\mathcal E))$ given by
\begin{equation}
\langle s_1, s_2 \rangle=\int _M g( s_1, s_2)\nu _g,\ \ \ \langle
S_1, S_2 \rangle=\int _M g( S_1, S_2 )\nu_g
\end{equation}
for sections $s_1, s_2 $ of  $\mathcal E$ and  $S_1, S_2$ --
sections of $HOM(TM,\mathcal E)$. Let $\nabla ^*$ denote
 the operator adjoint to $\nabla$ relative to $\langle,\rangle$, that is,
$$\langle \nabla s,
S\rangle= \langle s, \nabla ^*S\rangle $$
 for each section $s$
of  $\mathcal E$ and each section $S$ of $HOM(TM, \mathcal E)$. By
the formula (\ref{iloczyn_skalarny_nabla_sS}) and the divergence
theorem we obtain
\begin{equation}
\int _M g( \nabla s, S)\nu _g =\int _M\left( -\sum _{i=1}^n
g(s,(\onabla _{e_i}S)e_i)+g(s,SE)\right)\nu _g.
\end{equation}
 It is
now justified to set
\begin{equation}
\nabla ^* S=-\tr _g(\onabla _{\cdot}S)(\cdot) +SE.
\end{equation}
The definition makes sense also in the case where $M$ is neither
compact nor oriented. If $S=\nabla s$, where $s\in \mathcal
C^{\infty}(M\leftarrow \mathcal E)$ then we have
\begin{equation}\label{nabla*_nabla}
\nabla^*\nabla s=-\tr _g(\onabla_{\cdot}(\nabla
s))(\cdot)+\nabla_Es.
\end{equation}
Observe that if for a fixed point $x_0$ we take a local orthonormal
frame $e_1,...,e_n$ around $x_0$ such that $\hnabla e_{i}=0$ at
$x_0$ then at $x_0$ we have
\begin{equation}\label{formula_nabla(e_i,e_i)}
\nabla ^*\nabla s=-\tr _g\onabla _\cdot(\nabla _\cdot s).
\end{equation}
Assume now that $(g, \nabla, \nu)$ is an equiaffine statistical
structure. We can consider a scalar product on  tensor fields
determined by a volume form $\nu$.  Namely,  we set
\begin{equation}
\langle  s_1,s_2\rangle _\nu=\int _Mg(s_1,s_2)\nu,\ \ \ \ \langle
S_1,S_2 \rangle _\nu =\int _M g(S_1,S_2)\nu
\end{equation}
for sections $s_1,s_2$ of $\mathcal E$ and $S_1,S_2$ -- sections of
$HOM(TM,\mathcal E)$. Denote by $\nabla ^{*\nu}$ the operator
adjoint to $\nabla$ relative to the scalar product $\langle,\rangle
_\nu$, that is,
$$\langle \nabla s , S\rangle _\nu =\langle s,\nabla
^{*\nu}S\rangle _\nu .$$ By (\ref{iloczyn_skalarny_nabla_sS_nabla})
and the divergence theorem for equiaffine structures we get
\begin{equation}\label{nabla*nu_nabla}
\nabla ^{*\nu}S=-\tr _g(\onabla_{\cdot} S)(\cdot)
\end{equation}
Therefore

\begin{equation}\label{relation}
\nabla^*\nabla s = \nabla ^{*\nu}\nabla s+\nabla _Es.
\end{equation}

\bigskip
\section{ Hodge-type  theorems for statistical structures}

  As in the case of Riemannian manifolds we have  a Hodge-type decomposition
theorem for  compact equiaffine statistical manifolds. Let  $M$ be a
compact  manifold endowed with an equiaffine statistical structure
$(g,\nabla ,\nu)$. In order to prove a decomposition theorem observe
first that we can use the classical Fredholm alternative for the
product $\langle ,\rangle _\nu$. Namely, we take the bundle
$\mathcal E =\Lambda ^k TM^* $ with the Euclidean metric on each
$\mathcal E _x$ given by $g_\nu (\omega,\eta)=g(\psi \omega,
\psi\eta)$, where $\nu =\psi^2\nu _g$. Then $g(\omega,\eta )\nu
=g_\nu (\omega,\eta)\nu_g$ and
$$\langle \omega,\eta \rangle _\nu
=\int _M g_\nu (\omega,\eta)\nu _g.$$ By (\ref{relation_laplacians})
it is clear that the  smooth linear differential operator $\Delta
^\nabla$ of order 2 is elliptic. By Corollary
\ref{Delta_nabla_self_conjugate} it is self-adjoint relative to
$\langle ,\rangle _\nu$. Therefore, by the Fredholm alternative,
(see suitable formulation, for instance, in \cite{J}), we know that
if $\eta$ is $\langle ,\rangle _\nu$-orthogonal to $\ker \Delta
^\nabla_{|\mathcal F^k(M)}$ then there is a unique smooth
differential form $\omega$ such that $\Delta ^\nabla \omega=\eta$.

Denote by $\mathcal H^{k,\nabla}(M)$ the space of all smooth
$\nabla$-harmonic  forms of degree $k$  on $M$. The infinite
dimensional space $\mathcal F^k(M)$ is equipped with the Euclidean
scalar product $\langle ,\rangle _\nu$. In this section the
orthogonality will mean the orthogonality  with respect to this
product.

Although one can just say that the following  decomposition theorem
follows from the theory of  elliptic differential operators we give
a proof of this theorem for making  considerations of this paper
complete. The theorem says as follows
\begin{thm}\label{Hodge_decomposition}
Let $M$ be a compact manifold equipped with an equiaffine structure
$(g,\nabla,\nu)$. We have the following decomposition of the space
$\mathcal F^k(M)$ into the  mutually orthogonal (relative to
$\langle ,\rangle _\nu$) subspaces for every $k=0,1,...$
\begin{equation}
\mathcal F^k(M) =d (\mathcal F^{k-1}(M)) \oplus \mathcal
H^{k,\nabla}(M)\oplus \delta^{\onabla }(\mathcal F ^{k+1}(M)).
\end{equation}
\end{thm}
\proof
 We have
the mapping
 \begin{equation}\label{mapping_h}
 h:\mathcal H^{k,\nabla}(M)\ni \omega
\to [\omega]\in H^k_{DR}(M),\end{equation}
 where $H^k_{DR}(M)$
stands for the $k$-th de Rham cohomology group of $M$. Let $\omega$
be $\nabla$-harmonic and exact, i.e. $\omega=d\eta$. Since the
operators of the exterior differential and $\onabla$-codifferential
are adjoint and $\omega$ is $\onabla$-coclosed we obtain
$$\langle \omega,\omega\rangle _\nu=\langle d\eta,\omega\rangle
_\nu=\langle \eta,\delta ^{\onabla} \omega\rangle _\nu =0.$$
Therefore $\omega=0$ and, consequently,   the mapping $h$ is
injective. It follows, in particular,  that the space $\mathcal
H^{k,\nabla}(M)$ is finite dimensional. Let $\beta _1,...,\beta_r$
be an orthonormal basis of this space. Define the mapping
$$H:\mathcal F^k(M)\to \mathcal H^{k,\nabla}(M)$$
given by $H\omega=\sum _{i=1}^r\langle\omega ,\beta_i \rangle
_\nu\beta_i$. For every $\omega\in \mathcal F ^{k}(M)$ we have
$\omega= H\omega + (\omega -H\omega)$ and it is easily seen that the
form $\omega-H\omega$ belongs to  $(\mathcal H^{k,\nabla}(M))^\bot$.
Hence
$$\mathcal F^k(M)=\mathcal H^{k,\nabla}(M) \oplus (\mathcal
H^{k,\nabla}(M))^\bot .$$ We now observe that the following
subspaces of $\mathcal F^k(M)$ are mutually orthogonal
$$\mathcal H^{k,\nabla}(M), \ \ \ \ \ \ d(\mathcal F^{k-1}(M)),\ \ \
\  \delta^{\onabla} (\mathcal F^{k+1}(M)).$$ Namely, let $\omega \in
\mathcal H ^{k,\nabla}(M)$, $\eta \in \mathcal F^{k-1}(M)$, $\mu
\in\mathcal F^{k+1}(M)$. We have the following  obvious equalities
$$\langle \omega ,d\eta\rangle _\nu=\langle \delta ^{\onabla}
\omega,\eta\rangle _\nu =0,\ \ \ \langle \omega,\delta^{\onabla}
\mu\rangle _\nu=\langle d\omega,\mu\rangle _\nu=0,\ \ \ \langle
d\eta, \delta ^{\onabla} \mu \rangle_\nu=\langle  d^2\eta,\mu
\rangle _\nu=0.$$ Thus we have the direct orthogonal  sum
$$d(\mathcal F^{k-1}(M))\oplus \mathcal H^{k,\nabla}(M) \oplus
\delta ^{\onabla} (\mathcal F ^{k+1}(M))\subset\mathcal F^k(M).$$

Finally we observe that for every form $\omega '\in (\mathcal
H^{k,\nabla}(M))^\bot$ the form belongs to $\d(\mathcal
F^{k-1}(M))\oplus \delta^{\onabla }(\mathcal F^{k+1}(M))$. Since
$\omega ' $ is orthogonal to $\ker \Delta^\nabla _{|\mathcal F^k(M)}
=\mathcal H^{k,\nabla}(M)$, by the remarks made before this theorem
we know that there is a  smooth $k$-form $\eta$ such that $\Delta
^\nabla \eta =\omega '$. We now have
$$\omega ' =d(\delta^{\onabla}\eta) +\delta^{\onabla} (d\eta) \in
d(\mathcal F^{k-1}(M))\oplus \delta^{\onabla }(\mathcal
F^{k-1}(M)),$$ which finishes the proof. \koniec

Observe now that the mapping $h$ given by (\ref{mapping_h}) is also
surjective. Indeed, let $[\omega]\in H^k_{DR}(M)$. There exist
$\eta\in \mathcal F^{k-1}(M)$, $\omega_H\in \mathcal H^{k,\nabla
}(M)$ and $\mu \in \mathcal F^{k+1}(M)$ such that
$$\omega = d\eta +\omega_H +\delta^{\onabla} \mu.$$
Since $0=d\omega =d\delta ^{\onabla} \mu$, we have $0=\langle
d\delta^{\onabla} \mu,\mu\rangle _\nu =\langle \delta^{\onabla}\mu,
\delta^{\onabla} \mu \rangle_\nu$. Hence $\delta ^{\onabla} \mu=0$
and consequently $\omega =d\eta+\omega _H$, that is,
$[\omega]=[\omega_H]$. Hence the following  representation theorem
holds
\begin{thm}\label{Hodge-representation}
Let $M$ be a compact manifold equipped with an equiaffine
statistical structure $(g,\nabla, \nu)$.  The mapping $h$ given by
{\rm (\ref{mapping_h})} is an isomorphism. In particular, $\dim
\mathcal H^{k,\nabla}(M)=b_k(M)$, where $b_k(M)$ is the $k$-th Betti
number of $M$.
\end{thm}

\bigskip
\section{Bochner's technique for vector fields and harmonic $1$-forms.}

We shall first collect  basic information concerning vector fields
and their dual $1$-forms  on statistical manifolds. In this section
we assume that $(g,\nabla)$ is a statistical structure on $M$. As in
Section 2 we set $S_X=\nabla X$. Analogously $\hat S_X=\hat\nabla X$
and $\overline S_X=\onabla X$.
\begin{lemma}
Let $(g,\nabla)$ be a statistical structure on $M$. For a vector
field $X\in \mathcal X(M)$ the $(1,1)$-tensor field $S_X$ is
symmetric relative to $g$ if and only if $dX^{\flat}=0$.
\end{lemma}
\proof Let  $\eta =X^\flat$.  Using also the assumption that
$\nabla$ is torsion-free, we get
\begin{eqnarray*}
 d\eta (U,V)&=& U(\eta (V))-V(\eta(U)) -\eta([U,V])\\
&=& U(g(V,X))- V(g(U,X)) -g(\nabla _UV,X) +g(\nabla _VU, X)\\
&=& (\nabla _Ug) (V,X) +g(V,\nabla _UX)-(\nabla _Vg)(U,X) -g(U,\nabla _VX)\\
&=& g(V,S_XU)-g(U,S_XV).
\end{eqnarray*}\koniec

Let $X\in \mathcal X(M)$. Since $\delta X^\flat =-\tr _g\hnabla
X^\flat (\cdot,\cdot)=-\tr \hat S _X$, we have $-\tr S_X =-\tr \hat
S_X-\tau(X)=\delta X^\flat-\tau (X).$
 Since $X^\flat(E)=\tau (X)$ and $\delta ^{\onabla}=\delta -\iota_E$, one gets

\begin{lemma}\label{divX}
Let $(g,\nabla)$ be a statistical structure on $M$. For any $X\in
\mathcal X(M)$ we have
\begin{equation}
-\tr S_X=\delta X^\flat -\iota _E X^{\flat}=\delta
^{\onabla}X^\flat.
\end{equation}
\end{lemma}



\begin{proposition}\label{closed_coclosed_1-forms}
Let $(g,\nabla)$ be a statistical structure on $M$ and $X\in
\mathcal X(M)$. The $1$-form $X^{\flat}$ is closed if and only if
$S_X$ is symmetric relative to $g$. The $1$-form $X^{\flat}$ is
coclosed if and only if $\tr S_X=\tau(X)$. The $1$-form is
$\onabla$-coclosed if and only if $\tr S_X=0$.
\end{proposition}

\begin{lemma}\label{nablaX=0_onablaX=0} Let $(g,\nabla)$ be a statistical structure on $M$. For any
  $X\in \mathcal X(M)$ we have: $g( \nabla X, \nabla
X)=g(\onabla X ^\flat,\onabla X ^\flat)$. In particular, $\nabla
X=0$ if and only if $\onabla X ^\flat =0$.
\end{lemma}
\proof Let  $\eta =X^\flat$. It is sufficient to observe that
\begin{equation}\label{onabla_eta_nablaX}
\onabla \eta (U,V)=U(\eta(V))-\eta (\onabla_
UV)=U(g(X,V))-g(X,\onabla_UV)=g(\nabla _UX,V).
\end{equation}\koniec

By duality $\nabla \eta=0$ if and only if $\onabla \eta ^\sharp =0$
for any $1$-form $\eta$.
\begin{corollary}
If $(g,\nabla)$ is a trace-free statistical structure on $M$ and
$\nabla X=0$ for $X\in \mathcal X(M)$ then $X^{\flat}$ is closed and
coclosed. If for a $1$-form $\eta$ we have $\nabla\eta=0$ then
$\eta$ is closed and coclosed.
\end{corollary}

In some situations harmonic forms are  parallel. A well-known
Bochner theorem says that  harmonic 1-forms are parallel relative to
the Levi-Civita connection  on a Ricci non-negative Riemannian
manifold. We shall now prove some generalizations of this theorem.

\begin{thm}\label{thm_nabla_harmonic_1forms}
Let $M$  be a connected compact oriented manifold with an equiaffine
statistical structure
 $(g,\nabla, \nu)$.
If the Ricci tensor $Ric$ for $\nabla$ is non-negative on $M$ then
every $\nabla$-harmonic $1$-form on $M$ is $\onabla$-parallel. In
particular, the first Betti number $b_1(M)$ is not greater than
$\dim M$. If additionally $Ric>0$ at some point of $M$ then
$b_1(M)=0$.
\end{thm}
\proof Let $\eta$ be a $\nabla$- harmonic $1$-form and $X=\eta
^{\sharp}$. Since $M$ is compact,  $\eta$ is closed and
$\onabla$-coclosed.  Proposition \ref{closed_coclosed_1-forms}
yields that $\tr S_X=0$ and $S_X$ is symmetric relative to $g$.  In
particular, $S_X$ is diagonalizable. Therefore $\tr S_X^2\ge 0$ and
the equality holds if and only if $S_X=0$. By Lemma
\ref{th_int_Ric(X,Y)} we have
$$\int _MRic(X,X)\nu =-\int _M\tr S_X^2\nu\le 0.$$ If $Ric\ge 0$ on $M$ then
 $\tr S_X^2=0$ on $M$  and consequently $S_X$ vanishes on $M$. By
 Lemma \ref{nablaX=0_onablaX=0} we get
$\onabla \eta =0$.

 Let $x$ be any point of $M$ and
consider the mapping sending each  $\nabla$-harmonic $1$-form $\eta$
to $\eta _x\in T_x M^*$. The mapping is linear and, since each
$\nabla$-harmonic $1$-form is covariant constant, the mapping is
also a monomorphism. By Theorem \ref{Hodge-representation}
$b_1(M)\le \dim M$.

 We also  have $Ric (X,X)=0$ on $M$. 
 Hence, if $Ric>0$ at some point then $X=0$ at this point
and consequently, since $X$ is covariant constant,  $X=0$ on $M$.
Consequently $\eta =0$ on $M$.\koniec

In particular, if $(g,\nabla, \gnu)$ is equiaffine we get

\begin{corollary}\label{Bochner_harmonic}
Let $M$ be a connected compact oriented manifold with a trace-free statistical structure
 $(g,\nabla)$.
If the Ricci tensor $Ric$ for $\nabla$ is non-negative on $M$ then
every harmonic $1$-form on $M$ is $\onabla$-parallel. In particular,
the first Betti number $b_1(M)$ is not greater than $\dim M$. If
additionally $Ric>0$ at some point of $M$ then $b_1(M)=0$.
\end{corollary}



If $(g,\nabla)$ is a trace-free statistical structure and  $X^\flat$
is harmonic then $\tr S_X=0$ and $\tr S_X^2=g(\nabla X, \nabla X)$.
Therefore, if $M$ is compact and oriented, by Lemma
\ref{th_int_Ric(X,Y)} we have
\begin{equation}
\int Ric(X,X) \gnu =-\int g(\nabla X,\nabla X)\gnu .
\end{equation}
The same formula holds for the conjugate connection $\onabla$.

In the same way as Theorem \ref{thm_nabla_harmonic_1forms} one gets
\begin{thm}\label{glupie_uogolnienie}
Let M be a connected compact oriented manifold. Let $(g,\nabla)$  be
a trace-free statistical structure on $M$. If $Ric +\overline{Ric}
\ge 0$ on $M$ then each harmonic 1-form on $M$ is parallel relative
to the connections $\nabla$, $\onabla$ and $\hnabla$. In particular,
$b_1(M)\le \dim M$. If moreover $Ric+\overline{Ric}>0$ at some
point then $b_1(M)=0$.
\end{thm}

Note that the assumption  $Ric+\overline{Ric}\ge 0$ implies
that $\widehat{Ric}\ge 0$ and the classical Bochner theorem implies that
 harmonic $1$-forms are $\hnabla$-parallel but it does not imply
that they are $\nabla$ or $\onabla$-parallel.

\medskip

Consider now the case where $(g,\nabla, \nu )$ is equiaffine and a
vector field $X$ equals to $X=grad\, f$ for some function $f\in
\mathcal C ^\infty (M)$.  Let $\eta =df$. We have
\begin{equation}\label{onabla_eta=Hess}
\onabla \eta =\onabla ^2 f=Hess ^{\onabla}f.\end{equation}
 Since $Hess ^{\onabla}$ is a
symmetric $(0,2)$-tensor field (consequently diagonalizable at each
point of $M$)  and $\Delta ^\nabla f=-\tr _g Hess^{\onabla}f$, by
the Schwarz inequality we  have
\begin{equation}\label{Schwarz_inequality}
n|Hess^{\onabla}f|^2\ge|\Delta^{\nabla} f|^2.
\end{equation}
Since $df$ is closed, $S_X$ is symmetric (by Proposition
\ref{closed_coclosed_1-forms}) and therefore $\tr S_X
^2=g(S_X,S_X)=g(\nabla X,\nabla X)$. By Lemma
\ref{nablaX=0_onablaX=0} and the formula (\ref{onabla_eta=Hess}) we
now get
\begin{equation}
\tr S_X ^2=|Hess ^{\onabla} f|^2.
\end{equation}

 We can now prove the following generalization
 of a Bochner-Lichnerowicz formula and Lichnerowicz's theorem

\begin{thm}
Let $M$ be a compact  manifold. If $(g,\nabla,\nu)$ is an equiaffine
statistical structure on $M$ then for every function $f\in\mathcal C
^{\infty}(M)$ we have

 \begin{equation}\label{nabla-formula_Lichnerowicz}
 \int _M Ric (df^\sharp, df^\sharp)\nu= \int _M|\Delta ^\nabla f|^2 \nu-\int
 _M|Hess ^{\onabla} f|^2\nu
 \end{equation}
 for any function $f\in\mathcal C^\infty (M)$.
If for some real number $k$
$$ Ric \ge kg,$$
then the first eigenvalue $\lambda _1$ of the Laplacian $\Delta
^\nabla$ satisfies the inequality
\begin{equation}\label{inequality} \lambda _1\ge
{n\over{n-1}}k.\end{equation}\end{thm}

\proof The equality (\ref{nabla-formula_Lichnerowicz}) immediately
follows from (\ref{basic}) and the fact that $\tr S_X =-\Delta ^
\nabla f$ (by (\ref{onabla_eta_nablaX})), where $X=df^\sharp$. If
$\lambda$ is an eigenvalue of the Laplacian $\Delta ^\nabla$ and $f$
is the corresponding eigenfunction then (by Lemma
\ref{d_delta_onabla_adjoint})
$$\int_M g(\Delta ^\nabla f,\Delta ^\nabla
f)\nu =\int_M\lambda g(f,\Delta ^\nabla f)\nu=\int_M\lambda
g(f,\delta ^{\onabla}df)\nu=\int_M\lambda g(df, df)\nu.$$
 The second statement now  follows from
the formulas ({\ref{nabla-formula_Lichnerowicz}}),
(\ref{Schwarz_inequality}) and the assumed inequality. \koniec

The above theorem was proved in \cite{W} by using a different
method.
\begin{corollary}
Let $M$ be a compact oriented  manifold. If $(g,\nabla)$ is a
trace-free  statistical structure on $M$ then for every function
$f\in\mathcal C^\infty (M)$ we have
 \begin{equation}\label{formula_Lichnerowicz}
 \int _M Ric (df^\sharp, df^\sharp)\gnu= \int _M|\Delta  f|^2\gnu -\int
 _M|Hess ^{\onabla} f|^2\gnu
 \end{equation}
 for any function $f\in\mathcal C ^\infty (M)$.
If for some real number $k$
$$Ric \ge kg,$$
then the first eigenvalue $\lambda _1$ of the Laplacian $\Delta $
satisfies the inequality
$$ \lambda _1\ge {n\over{n-1}}k.$$
\end{corollary}

\begin{remark}{\rm In the Riemannian case it is known that the equality in (\ref{inequality}) holds
if and only if $(M,g)$ is isometric to an ordinary sphere. In the
case of statistical structures we have the following
\begin{example}
{\rm Let $\mathbf{f}:M\to \R^{n+1}$ be a locally strongly convex
immersed hypersurface. Equipping it with the transversal vector
field $\xi =-\mathbf{f}$ we get the induced statistical  structure
on $M$. We say that $\mathcal f:M\to \R^{n+1}$ is  a centroaffine
hypersurface. If $\alpha$ is a $1$-form on $\R ^{n+1}$ we define the
function $f$ on $M$ by the formula $f=\alpha (\xi)=
-\alpha\circ\mathbf{f}$. Denote by $\nabla$ the induced connection.
We have

\begin{eqnarray*}
(\nabla df) (X,Y) &=&X(df(Y)) -df(\nabla_XY) =-\alpha (X(Y\mathbf
{f}))+\alpha (\mathbf{f}_* (\nabla _XY))\\
 &=& -\alpha
(\mathbf{f}_*(\nabla_XY) -g(X,Y)\mathbf{f})
+\alpha(\mathbf{f}_*(\nabla_XY))\\&=&-f g(X,Y).
\end{eqnarray*}

Hence $\nabla df+fg=0$. Apply this to the case where  $\mathbf {f}$
is the conormal map of some centroaffine immersion and $g$ is the
common second fundamental form for the immersion and its conormal.
Then for functions $f$ defined as above we have
\begin{equation}\label{Obata} \onabla df +fg =0,\end{equation}
that is, $\Delta ^\nabla f=f$. We propose the following conjecture:
 {\it     Let $M$ be a compact manifold equipped with a statistical structure $(g,\nabla)$.
 Assume that the conjugate connection $\onabla$ is complete.
 If there is a non-constant  function $f$ such that the equation {\rm(\ref{Obata})}
 is satisfied then $(g,\nabla )$ can be realized as the induced structure  on a centroaffine ovaloid.}
}\end{example} }
\end{remark}

We continue considerations of the section assuming that  a
statistical structure $(g,\nabla)$ is trace-free.
Recall the classical Bochner formula for the Laplacian of the square
of the length of a vector field $X$ for which $dX ^\flat=0$
(equivalently $S_X$, $\overline S_X$, $\hat S_X $ are symmetric
relative to $g$). Namely, if $\varphi= g(X,X)$ then
\begin{equation}\label{Bochner_hnabla}
\Delta\varphi + 2X(\div^{\gnu} X)=-2\widehat{ Ric}(X,X)- 2g(\hnabla
X,\hnabla X).
\end{equation}
On the other hand, by formula (\ref{lematRic(X,Y)}) applied for
$\hat\nabla$
 we have
$$2\div^ {\gnu}(\hnabla _XX)=2\widehat{Ric} (X,X)+2X(\div^{\gnu} X) +2g(\hnabla X,
\hnabla X).$$
Adding these two equalities we get
\begin{equation}
\Delta \varphi =-2\div^{\gnu} \hnabla_XX.
\end{equation}
Since $S_X$ is symmetric, $\tr S_X ^2=g(\nabla X,\nabla X)$.
Similarly $\tr \overline S_X^2= g(\onabla X,\onabla X)$.  By
(\ref{lematRic(X,Y)}) applied for $\nabla$, $\onabla$ and $\hnabla$
and the fact that $\tr S_X=\tr \hat S_X=\tr \overline
S_X=\div^{\gnu} X$ we get
$$ \div^{\gnu}(\nabla _XX)=Ric(X,X) +Xdiv ^{\gnu}X + g(\nabla X,\nabla X),$$
$$ \div^{\gnu} (\onabla _XX)=\overline{Ric}(X,X) +Xdiv ^{\gnu}X + g(\onabla X,\onabla X).$$
Since $2\hnabla _XX=\nabla _XX+\onabla _XX$, we have
\begin{equation}\label{Bochner_nabla_onabla}
\Delta \varphi +2Xdiv^{\gnu} X = -Ric(X,X)-\overline{Ric} (X,X) -
g(\nabla X,\nabla X) -g(\onabla X,\onabla X)
\end{equation}
for any vector field $X$ such that $dX^\flat =0$. By
(\ref{g(nabla,nabla)+g(onabla,onabla)}) the last formula can be
equivalently written as
\begin{equation}\label{Bochner_hnabla_K}
\Delta \varphi +2Xdiv^{\gnu} \, X = -Ric(X,X)-\overline{Ric} (X,X) -
2g(\hnabla X,\hnabla X) -2g(K_X,K _X).
\end{equation}

\begin{thm}
Let $(g,\nabla)$ be a trace-free statistical structure on a connected manifold  $M$
 and $Ric+\overline{Ric}\ge 0$ on $M$.
Let $\eta$ be a  closed  harmonic  $1$-form on $M$.
\newline
{\rm 1)} If $\varphi =|\eta|^2$ attains a local maximum at some
point $x_o$  of $M$ then $\nabla\eta=\onabla \eta=\hnabla\eta=0$ at
$x_o$. If moreover $Ric +\overline Ric >0$ at $x_0$ then $\eta=0$ in
a neighborhood of $x_o$.
\newline
{\rm 2)} If $Ric+\overline{Ric}>0$ on $M$ and $\varphi$ attains a
global maximum at some point of $M$ then $\eta =0$ on $M$.
\newline
{\rm 3)} If $Ric+\overline{Ric}>0$ on $M$, $\varphi$ attains a local
maximum at some point and  $g$ is analytic then $\eta=0$ on $M$.
\end{thm}
\proof Let $X=\eta^\sharp$. Of course $\varphi=g(X,X)$. If $\varphi$
attains a local maximum at $x_o$ then $(\Delta\varphi)_{x_o}\ge 0$.
Since $\eta $ is closed and harmonic, we have that  $\delta \eta$ is
constant. By Lemma \ref{divX} it follows that $X (\div^{\gnu} X)=0$.
Using now (\ref{Bochner_nabla_onabla}) we  obtain the first
assertion in {\rm 1)}. Moreover $Ric (X,X) +\overline {Ric}(X,X) =0$
at $x_0$. If $Ric+ \overline {Ric}>0$ at $x_0$ then $X_{x_o}=0$,
i.e., $\eta_{x_o}=0$. Since $\varphi$ attains a local maximum at
$x_o$, we have that $\eta=0$ around $x_o$. If the maximum is global
then $\eta$ vanishes on $M$. Since a harmonic form on an analytic
Riemannian manifold is analytic, we have {\rm 3)}. \koniec

\section{Bochner-Weitzenb\"ock formulas for differential forms}
Let $(g,\nabla)$ be a statistical structure on a manifold $M$. The
Weitzenb\"ock  curvature operator for the curvature tensor $R$ for
$\nabla$ will be denoted by $\mathcal W ^R$ . Let $s$ be a tensor
field of type $(l,k)$, where $k>0$, on $M$. One defines a tensor
field $\mathcal W^Rs$ of type $(l,k)$ as follows

\begin{equation}
(\mathcal W^Rs)(X_1,...,X_k)= \sum
_{i=1}^k\sum_{j=1}^n(R(e_j,X_i)s)(X_1,...,e_j,..., X_k),
\end{equation}
where $e_1,...,e_n$ is an arbitrary orthonormal frame, $R(e_j,X_i)s$
means that $R(e_j,X_i)$ acts as a differentiation on $s$, and $e_j$
in the last parenthesis is at the $i$-th place. The definition is
 independent of the choice of  an orthonormal basis.

 Observe
what the $\mathcal W^R\omega$ is in the case where $\omega$ is a
$1$-form. In this case we have
\begin{eqnarray*}
\sum _{j=1}^n(R(e_j,X)\omega)(e_j)&=&-\sum
_{j=1}^n\omega(R(e_j,X)e_j)=-\sum _{j=1}^n g(R(e_j,X)e_j,\omega^\sharp)\\
&=& \sum _{j=1}^n g(\overline R (e_j,X)\omega^\sharp,
e_j)=\overline{Ric}(X,\omega^\sharp).
\end{eqnarray*}
Thus for $1$-forms
\begin{equation}\label{Weitzenbock_for_1-forms}
\mathcal W^R\omega(X)=\overline{Ric}(X,\omega ^{\sharp}).
\end{equation}

We shall now prove some generalizations of the Bochner-Weitzenb\"ock
formula for the Laplacians acting on differential forms on
statistical manifolds.

\begin{thm}\label{main_theorem}{\rm i)} For any statistical structure  $(g,\nabla)$ we have
\begin{equation}\label{i}
\Delta=\nabla^*\nabla +\mathcal W^R +\nabla _E -\mathcal L_E=
\nabla^*\nabla +\mathcal W^R +S_E,
\end{equation}
\begin{equation}\label{i_for_Delta_nabla}
\Delta^\nabla=\nabla^*\nabla +\mathcal W ^R +2S_E-\nabla _E.
\end{equation}
\newline
\noindent {\rm ii)} If $(g,\nabla, \nu)$ is an equiaffine
statistical structure then
\begin{equation}\label{ii}
\Delta ^\nabla=\nabla^{*\nu}\nabla +\mathcal W^R +2\nabla _E
-2\mathcal L_E=\nabla^{*\nu}\nabla +\mathcal W^R +2S_E.
\end{equation}
\newline
\noindent {\rm iii)} If $(g,\nabla)$ is  trace-free then
\begin{equation}\label{iii}
\Delta =\nabla^*\nabla +\mathcal W^R.
\end{equation}

\end{thm}
\proof i) Let $\omega$ be a $k$-form on $M$. Let $x_o$ be a fixed
point of $M$ and $e_1,...,e_n$ be a local orthonormal frame around
$x_o$ such that $\hat\nabla e_i=0$ at $x_o$. As before, we shall use
the fact that  $\sum_{i=1}^n\nabla _{e_i} e_i=E$ at $x_o$. Let
$X_1,...,X_k$ be arbitrary vectors from $T_{x_o}M$. We extend them
to local vector fields around $x_o$  in such a way that $\hat \nabla
X_i=0$ at $x_o$ for $i=1,...,k$. Then $(\nabla
_YX_i)_{x_o}=(K_YX_i)_{x_0}$ and $[e_j,X_i]_{x_o}=0$ for each $Y$
and every $i=1,...,k$ and $j=1,...,n$.


 We shall now compute at
$x_o$
\begin{eqnarray*}
d\delta ^\nabla\omega(X_1,...,X_k)&=&\sum_{i=1}^k(-1)^{i-1}
(\nabla_{X_i}(\delta^\nabla\omega))(X_1,...,\hat{X_i},...,X_k)\\
&=&\sum _{i=1}^k(-1)^{i-1}\{X_i((\delta^\nabla\omega)(X_1,...,\hat
X_i,...,X_k))\\
&&-\delta^\nabla \omega(\nabla _{X_i}X_1,...,\hat {X_i},...,X_k)\\
&&\ \ \ \  \ \-.......-\delta^\nabla\omega(X_1,...,\hat
X_i,...,\nabla _{X_i}X_k)\}\\
&=& \sum_{i=1}^k\sum_{j=1}^n (-1)^iX_i((\nabla
_{e_j}\omega)(e_j,X_1,...,\hat X_i,...X_k))\\
&&+\sum_{j=1}^n\sum_{i=1}^k(-1)^{i-1}(\nabla_{e_j}\omega)(e_j,\nabla_{X_i}X_1,...,\hat
X_i,...,X_k)\\
&&\ \ \   +........
+\sum_{j=1}^n\sum_{i=1}^k(-1)^{i-1}(\nabla
_{e_j}\omega)(e_j,...,\hat X_i,...,\nabla _{X_i}X_k)\\
&=&\sum_{i=1}^k\sum_{j=1}^n(-1)^i(\nabla
_{X_i}(\nabla_{e_j}\omega))(e_j,X_1,...,\hat X_i,...,X_k)\\
&& +\sum_{i=1}^k\sum_{j=1}^n(-1)^i(\nabla _{e_j}\omega)(\nabla
_{X_i}e_j,X_1,...,\hat X_i,...,X_k)\\
&&\  + \sum_{i=1}^k\sum_{j=1}^n(-1)^i(\nabla
_{e_j}\omega)(e_j,\nabla _{X_i}X_1,...,\hat X_i,...,X_k)\\
&&\ \ \ \ \ \
+.......+\sum_{i=1}^k\sum_{j=1}^n(-1)^i(\nabla_{e_j}\omega)(e_j,...,\hat
X_i,...,\nabla_{X_i}X_k)\\
&&+\sum_{i=1}^k\sum_{j=1}^n(-1)^{i-1} (\nabla
_{e_j}\omega)(e_j,\nabla_{X_i}X_1,...,\hat X_i,..., X_k)\\
&&\ \ \ \ \ \   +.......+
\sum_{i=1}^k\sum_{j=1}^n(-1)^{i-1}
(\nabla _{e_j}\omega)(e_j,...,\hat X_i,...,\nabla_{X_i}X_k)\\
&=&-\sum_{i=1}^k\sum_{j=1}^n(\nabla _{X_i}(\nabla
_{e_j}\omega))(X_1,...,e_j,...,X_k)\\
&&-\sum_{i=1}^k\sum_{j=1}^n(\nabla _{e_j}\omega)(X_1,...,K
_{e_j}X_i,...,X_k),
\end{eqnarray*}
where in the last parenthesis in the last two lines $e_j$ and $K
_{e_j}X_i$  are at the $i$-th place. Using Lemma \ref{in_the_long_proof} for $\alpha =d\omega$ and
$\alpha =(\nabla _{X_i}\omega)(\cdot,...\hat
X_i,..., \cdot)$ (if $\omega$ is a $1$-form some lines in the following formula do not appear)
  we continue computations at
$x_o$
\begin{eqnarray*}
\delta d\omega(X_1,...,X_k)&=&-\sum _{j=1}^n(\nabla
_{e_j}(d\omega))(e_j,X_1,...,X_k)-(\iota_Ed\omega)(X_1,...,X_k)\\
&&= -\sum _{j=1}^n e_j((d\omega)(e_j,X_1,...,X_k))\\
&&\ \ \ \ +\sum _{j=1}^n d\omega (\nabla _{e_j}e_j,X_1,...,X_k)\\
&&\ \ \ \ +\sum _{j=1}^n d\omega(e_j,\nabla _{e_j}X_1,...,X_k)\\
&& \ \ \  + ..........
+\sum _{j=1}^n d\omega(e_j,...,\nabla _{e_j}X_k)- (\iota_Ed\omega)(X_1,...,X_k)\\
&=&-\sum _{j=1}^ne_j((\nabla _{e_j}\omega)(X_1,...,X_k))\\
&&\ -\sum _{i=1}^k\sum
_{j=1}^n(-1)^ie_j((\nabla_{X_i}\omega)(e_j,X_1,...,\hat
X_i,...,X_k))\\
&=&-\sum _{j=1}^n(\nabla _{e_j}(\nabla_{e_j}\omega))(X_1,...,X_k)\\
&&\ \ \ \ \ \ -\sum _{j=1}^n(\nabla_{e_j}\omega)(\nabla
_{e_j}X_1,...,X_k)\\&&\ \ \ \ \ \ \ \ \ \ \ \ \ - ..........-\sum
_{j=1}^n(\nabla_{e_j}\omega)(X_1,...,\nabla _{e_j}X_k)\\
&&- \sum _{i=1}^k \sum _{j=1}^n(-1)^i\{(\nabla _{e_j} (\nabla
_{X_i}\omega))(e_j,X_1,...,\hat X_i,...,X_k) \\&&\ \ \ \ \ \ \ \ \ \
\ +(\nabla
_{X_i}\omega)(\nabla _{e_j}e_j,X_1,...,\hat X_i,...,X_k)\\
&&\ \ \ \ \ \ \  \ \ \ \ \ \ \ \ +(\nabla _{X_i}\omega)(e_j,\nabla
_{e_j}X_1,...,\hat X_i,...,X_k)\\ &&\ \ \ \ \ \ \ \ \ \ \ \ \ \ \ \
\  \ +........... +(\nabla
_{X_i}\omega)(e_j,X_1,..,\hat X_i,...,\nabla _{e_j}X_k)\}\\
&=& -\sum _{j=1}^n(\nabla _{e_j}((\nabla_{e_j}\omega))(X_1,...,X_k)\\
&&\ \ \ \
\  \ \ \ -\sum _{j=1}^n\sum_{i=1}^k(\nabla
_{e_j}\omega)(X_1,..,K _{e_j}
X_i,...,X_k)\\
&& \ \ \ \ \ \ \ \ \ \ \ \ +\sum _{i=1}^k\sum
_{j=1}^n(\nabla_{e_j}(\nabla _{X_i}\omega))(X_1,...,e_j,...,X_k)\\
&&\ \ \ \ \ \ \ \ \ \ \ \ \ \ \ \ \ +d(\iota_E\omega)(X_1,...,X_k) + (S_E\omega) (X_1,...,X_k),
\end{eqnarray*}
where in the last line $e_j$ is at the $i$-th place and in the last
but one line $K _{e_j}X_i$ is at the $i$-th place. Finally we
observe that at $x_0$ we have

\begin{eqnarray*}
\sum _{j=1}^n (\onabla _{e_j}(\nabla
_{e_j}\omega))(X_1,...,X_k)&=&\sum_{j=1}^n((\nabla _{e_j}-2K_{e_j})(\nabla
_{e_j}\omega))(X_1,...,X_k)\\
&&= \sum _{j=1}^n(\nabla _{e_j}(\nabla _{e_j}\omega))(X_1,...,X_k)\\
&&\ \ \ \ \ \ \ \ +2\sum_{j=1}^n\sum_{i=1}^k (\nabla
_{e_j}\omega)(X_1,...,K_{e_j}X_i,...,X_k),
\end{eqnarray*}
where in the last line $K_{e_j}X_i$ is at the $i$-th place.
Composing the last three  formulas and using
(\ref{formula_nabla(e_i,e_i)}) we obtain the equality
\begin{equation}\label{formula_with_S}
\Delta\omega=\nabla^*\nabla\omega +\mathcal W^R \omega +S_E \omega
\end{equation}
Using (\ref{about_L_S}) and (\ref{relation_laplacians}) completes
the proof of i). To prove ii) it is sufficient to use i),
(\ref{relation_laplacians}) and (\ref{relation}).
The last statement is a particular case of i).
 \koniec

Using  (\ref{Weitzenbock_for_1-forms}), the last theorem  and the
formula (\ref{nablag_gK}) we get
\begin{corollary}
Let $(g,\nabla) $ be a statistical structure. For any $1$-form
$\omega$ we have
\begin{equation}
\Delta \omega=\nabla ^*\nabla \omega +\overline{Ric}(\cdot,\omega
^{\sharp})+\nabla _E\omega-\mathcal L_E\omega ,
\end{equation}
\begin{equation}
(\Delta  \omega)(X)=-\tr _g(\nabla _{\cdot,\cdot}^2\omega )(X)
+\overline{Ric}(X,\omega ^\sharp)+ g(\nabla _{X}g, \nabla _E\omega)
+\nabla_E\omega-\mathcal L _E\omega
\end{equation}
for any vector $X\in T_xM$, $x\in M$. If $(g,\nabla,\nu)$ is an
equiaffine statistical structure then
\begin{equation}
\Delta ^\nabla \omega=\nabla ^{*\nu}\nabla \omega
+\overline{Ric}(\cdot,\omega ^{\sharp})+2\nabla _E\omega-2\mathcal
L_E \omega.
\end{equation}
\end{corollary}

The formula (\ref{i}) can also be written  for the connection
$\onabla$.  Namely, we have
\begin{equation}\label{i_for_onabla}
\Delta =\onabla ^*\onabla +\mathcal W ^{\overline R} +\overline
S_{\overline E}=\onabla ^*\onabla +\mathcal W ^{\overline R} - S_E
+2K_E.
\end{equation}
If $(g,\nabla,\nu )$ is an equiaffine statistical structure then, in
general, $(g,\onabla, \nu)$ is not equiaffine. It is equiaffine if
an only if $(g, \nabla)$ is trace-free.

The structure $(g,\onabla, \overline \nu)$ is equiaffine if and only
if $\overline \nu= e ^{-\rho}\nu$ where  $\tau = d\rho$. In
particular, $\tau$ must be exact.
In such a case we define
\begin{equation}\label{onabla*nu} \onabla^{*\overline \nu}S=-\tr _g(\nabla _{\cdot}S)  (\cdot)
\end{equation}
for a section $S$ of $HOM(TM,\mathcal E)$. If $M$ is compact then
$\onabla ^{*\overline\nu}$ is  the adjoint operator for $\onabla$
relative to $\langle, \rangle _{\overline \nu}$. For the same
reasons as (\ref{relation}) one gets
\begin{equation}
\onabla ^*\onabla=\onabla ^{*\overline\nu}\onabla -\onabla _E
\end{equation}
We can now compute

\begin{eqnarray*}
\Delta^\nabla=\Delta-\mathcal L_E&=&\onabla ^*\onabla+\mathcal W
^{\overline R}-\overline S_E-\mathcal L_E\\
 &=&\onabla ^{*\overline\nu}\onabla+\mathcal W ^{\overline
R}-\onabla_E-\overline S_E-\mathcal L_E\\
&=&\onabla ^{*\overline\nu}\onabla+\mathcal W ^{\overline
R}-2\onabla _E.\end{eqnarray*} We have proved
\begin{proposition}
If $(g,\nabla,\nu)$ is an equiaffine statistical structure and
$\tau$ is exact and equal to $d\rho$ then
\begin{equation}
\Delta^\nabla=\onabla ^{*\overline\nu}\onabla+\mathcal W ^{\overline
R}-2\onabla _E,
\end{equation}
where $\overline \nu =e^{-\rho}\nu$.
\end{proposition}

Consider again a  tensor vector bundle $\mathcal E$ over $M$. If $s$
is a section of $\mathcal E$ and $\tilde \nabla$ is a connection
then we have
\begin{eqnarray*}
\tilde\nabla^2_{Y,X}|s|^2&=&Y(d|s|^2(X))-d|s|^2(\tilde\nabla_YX)\\
&=& Y(X|s|^2)-d|s|^2(\tilde\nabla _YX)\\
&=& Y\{g(\nabla _Xs,s)+g(s,\onabla _Xs)\} -d|s|^2(\tilde\nabla_YX)\\
&=& g(\onabla _Y(\nabla_Xs),s) +g(\nabla_Xs,\nabla_Ys)\\
&&\ \ \ \ \ \ \ \ \ \ +g(\onabla_Ys,\onabla_Xs)
+g(s,\nabla_Y(\onabla_X s))-d|s|^2(\tilde\nabla_YX).
\end{eqnarray*}
Applying this formula to the connections $\hnabla$ and $\onabla$ one
gets

\begin{eqnarray*}
\hnabla ^2_{Y,X} |s|^2&=& g(\onabla _Y(\nabla_Xs),s)
+g(\nabla_Xs,\nabla_Ys)\\&&+g(\onabla_Ys,\onabla_Xs)
+g(s,\nabla_Y(\onabla_X s))-d|s|^2(\hnabla_YX)\end{eqnarray*} and

\begin{eqnarray*}
\onabla ^2_{Y,X} |s|^2&=& g(\onabla _Y(\nabla_Xs),s)
+g(\nabla_Xs,\nabla_Ys)\\&&+g(\onabla_Ys,\onabla_Xs)
+g(s,\nabla_Y(\onabla_X s))-d|s|^2(\onabla_YX).\end{eqnarray*}
  For an orthonormal frame $e_1,...,e_n$ such that $\hnabla
e_j=0$ at a fixed point $x_0$ we obtain  at this point
\begin{equation}\label{Delta_s2}
\Delta |s|^2 =-\sum _{i=1}^ng(\onabla_{e_i}(\nabla_{e_i} s),s)-
\sum_ {i=1}^ng( \nabla_{e_i}(\onabla_{e_i}s),s)-g(\nabla s,\nabla
s)-g(\onabla s,\onabla s)
\end{equation}
\begin{equation}\label{Delta_nabla_s2}
\begin{array}{rcl}
\Delta ^\nabla |s|^2 &=&-\sum _{i=1}^ng(\onabla_{e_i}(\nabla_{e_i}
s),s)- \sum_ {i=1}^ng( \nabla_{e_i}(\onabla_{e_i}s),s)\\
&&\ \ \ \ \ \  \ -g(\nabla s,\nabla s)-g(\onabla s,\onabla
s)-Eg(s,s).\end{array}
\end{equation}

Therefore, using (\ref{Delta_s2}), (\ref{Delta_nabla_s2}) and
(\ref{formula_nabla(e_i,e_i)}) for $\nabla$ and $\onabla$, we obtain

\begin{thm} For any statistical structure and any tensor field $s$
we have
\begin{equation}\label{Delta_s_2}
\Delta |s|^2 =g(\nabla ^*\nabla s,s)+ g( \onabla ^*\onabla
s,s)-g(\nabla s,\nabla s)-g(\onabla s,\onabla s)
\end{equation}
\begin{equation}\label{Delta_nabla_s_2}
\begin{array}{rcl}
\Delta^\nabla |s|^2& =&g(\nabla ^*\nabla s,s)+ g( \onabla ^*\onabla
s,s)\\ &&\ \ \ \ \ \ -g(\nabla s,\nabla s)-g(\onabla s,\onabla
s)-g(\nabla_Es,s)-g(\onabla_Es,s).\end{array}
\end{equation}
 In particular, if $\omega$ is a differential form then
\begin{equation}\label{Bochner_Weitzenbock_formula_for_s^2}
\Delta |\omega|^2 =2 g(\Delta \omega,\omega) -g(\mathcal W
^{R+\overline R}\, \omega, \omega) -g(\nabla\omega,\nabla\omega)
-g(\onabla\omega,\onabla\omega)-2g(K_E\omega,\omega).
\end{equation}
\begin{equation}\label{Bochner_Weitzenbock_for s^2_Delta_nabla}
\begin{array}{rcl}
\Delta ^\nabla|\omega|^2 &=&2 g(\Delta ^\nabla\omega,\omega)
-g(\mathcal W ^{R+\overline R}\, \omega, \omega)\\
&&\ \ \ \ \ \  -g(\nabla\omega,\nabla\omega)
-g(\onabla\omega,\onabla\omega)-2g(\nabla_E\omega,\omega).\end{array}
\end{equation}
\end{thm}
\proof To prove (\ref{Bochner_Weitzenbock_formula_for_s^2}) it is
sufficient to use (\ref{Delta_s_2}), (\ref{i}) and
(\ref{i_for_onabla}). Formula (\ref{Bochner_Weitzenbock_for
s^2_Delta_nabla}) follows from (\ref{Delta_nabla_s_2}) and
(\ref{i_for_Delta_nabla}).\koniec

From (\ref{Delta_s_2}) and the maximum principle applied to the
standard Laplacian we immediately get
\begin{proposition}
Let  $M$ be a connected manifold equipped with a  statistical
structure $(g,\nabla)$. If $s$ is a tensor field on  $M$ such that
\begin{equation}
g(\nabla ^*\nabla s+\onabla ^*\onabla s,s)\le 0
\end{equation}
on $M$ and $|s|$ attains a maximum then  $\nabla s=0$ and $\onabla
s=0$ on $M$.
\end{proposition}


\bigskip

\section{Bochner-Weitzenb\"ock formulas for the metric tensor field}
If $(g,\nabla)$ is a  statistical structure then we can  apply the
Weitzenb\"ock curvature operator $\mathcal W^R$
 to $g$. One easily sees that

\begin{equation}
(\mathcal W^Rg)(X,Y)=\overline {Ric}(X,Y)+ \overline{Ric }(Y,X)-Ric(X,Y)-Ric(Y,X).
\end{equation}
 For any
tensor field $s$ on $M$ we have $R(X,Y)s=\nabla ^2_{X,Y}s-\nabla^2
_{Y,X} s. $ Since $\nabla g$ is symmetric, we have
$(\nabla_{X,Y}^2g)(Z,W)=(\nabla _{X,Z}^2g)(Y,W)$. Therefore
$$\sum _{i=1}^n(R(e_i,X)g)(e_i,Y)=-\sum _{i=1}^n(\nabla _{X,e_i}^2g)
(e_i,Y)+\tr_g(\nabla_{\cdot,\cdot}^2g)(X,Y).$$

As usual we choose an orthonormal frame $e_1,...,e_n$ around a fixed
point $x_o$ such that $\hnabla e_i=0$ at $x_o$, that is, $\nabla _X
e_i=K_Xe_i$ at $x_o$ for any $X$. Using now formulas
(\ref{nablag_gK}) and  (\ref{tr_g_nablag_tau}) we get at $x_0$

\begin{eqnarray*} -\sum _{i=1}^n(\nabla
_{X,e_i}^2g)(e_i,Y)&=&-\sum_{i=1}^n(\nabla_X(\nabla
g))(e_i,e_i,Y)\\&=&\sum_{i=1}^n  2\nabla g(\nabla_X e_i,e_i,Y)
-X(\nabla  g(e_i,e_i,Y)) +\nabla g(e_i, e_i,\nabla _XY)\\
&=&  2\sum_{i,j=1}^n g(\nabla_X e_i,e_j) \nabla g(e_j,e_i,Y)+2\nabla\tau(X,Y)\\
&=&
2\sum_{i,j=1}^n g(K_X e_i,e_j) \nabla g(e_j,e_i,Y)+2\nabla\tau (X,Y)\\
&=& -\sum _{i,j=1}^n(\nabla _X g)(e_i,e_j) (\nabla
_Yg)(e_i,e_j)+2\nabla\tau(X,Y)\\
&=&-g(\nabla _Xg,\nabla_Yg)+2\nabla\tau (X,Y).
\end{eqnarray*}
On the other hand we have
\begin{eqnarray*}\sum _{i=1}^n(R(e_i,X)g)(e_i,Y)
=\overline {Ric}(X,Y)-Ric(X,Y).\end{eqnarray*}
 We have proved
\begin{lemma}
For any statistical connection $\nabla$ for $g$ we have
\begin{equation}\label{tr_gnabla_g^2}
\tr_g\nabla _{\cdot,\cdot}^2g(X,Y)-g(\nabla _Xg,\nabla
_Yg)+2\nabla\tau(X,Y)=-Ric(X,Y)+\overline{Ric}(X,Y),
\end{equation}
\begin{equation}
2\tr_g\nabla _{\cdot,\cdot}^2g(X,Y)-2g(\nabla _Xg,\nabla
_Yg)+2\nabla\tau(X,Y)+2\nabla\tau (Y,X)=(\mathcal W^Rg)(X,Y).
\end{equation}
\end{lemma}
Since the scalar curvature of $\nabla$ and the one for $\onabla$ are
the same we obtain
\begin{proposition}
For any trace-free statistical structure $(g, \nabla)$ we have $$
\sum_{i,j=1}^n\nabla^2g(e_i,e_i,e_j,e_j)=g(\nabla g,\nabla g).$$
\end{proposition}
An example of usage of this theorem is the following
\begin{corollary}
If for a locally strongly convex Blaschke hypersurface  for its
Blaschke metric $g$ and its induced connection $\nabla$
$$\sum_{i,j=1}^n\nabla^2g(e_i,e_i,e_j,e_j)=0,$$ then the hypersurface is a
quadric.
\end{corollary}
\proof It follows from Berwald's theorem and the above
proposition.\koniec

\bigskip

\section{The sectional curvature and the curvature operator for statistical structures}

\subsection{Algebraic preliminaries}

 Let $V$ be an
$n$-dimensional Euclidean vector space with the scalar product
$\langle \cdot ,\cdot \rangle$. Let $T$ be any $(1,3)$-tensor on
$V$. We also set $T(X,Y)Z:=T(X,Y,Z)$. If the tensor is
skew-symmetric relative to $X,Y$ and the  Bianchi identity:
$T(X,Y,Z)+T(Y,Z,X)+T(Z,X,Y)=0$ holds,  we  call $T$ a curvature-like
tensor of type $(1,3)$. For $T$ we define a $(0,4)$-tensor $ T$ as
follows
$$ T(U,Z,X,Y)=\langle T(X,Y)Z,U\rangle .$$
If $T$ is a curvature-like tensor of type $(1,3)$ and the tensor $
T$ of type $(0,4)$ is skew-symmetric relative to $U,Z$ then  the
both tensors $T$ will be  called  Riemann-curvature-like tensors.
For a Riemann-curvature-like tensor we have
$$T(U,Z,X,Y)=T(X,Y,U,Z).$$ The easiest Riemann-curvature-like tensor
is $R_0$ defined as follows
$$R_0(X,Y)Z=\langle Y,Z\rangle X-\langle X,Z\rangle Y.$$

The scalar product extended to the exterior products of $V$ gives the isometric
identification of $\Lambda ^2 V$
and $\Lambda ^2V^*$. The last space will be also isometrically
identified with the space  $\mathfrak{so}(V)$ of all skew-symmetric endomorphisms of  $V$.
In particular, if $T(X,Y)Z$ is skew-symmetric in $X,Y$, then $T(X\wedge Y)$ is well defined and consequently
$T(\Theta)$ is well-defined for any $\Theta \in \mathfrak{so}(V)$.
If $  T$ is a Riemann-curvature-like tensor then it defines the curvature  operator   $\mathfrak{T}$ sending
$2$-vectors  into $2$-vectors, that is,
\begin{equation}
\langle \mathfrak{T}(X\wedge Y), Z\wedge U\rangle=  T(X,Y,Z,U).
\end{equation}
Because of the  properties of $ T$,  the above formula
defines a linear, symmetric relative to  the given scalar
product, operator $\mathfrak{T}:\Lambda ^2 V  \to \Lambda ^2 V$.
In particular, $\mathfrak{T}$ is diagonalizable.

We shall now adapt the material of section 3 from \cite{P1} to the case we study.
Let $\Theta _{\alpha}$ be an orthonormal basis of $\Lambda ^2V$. For any Riemann-curvature-like tensor $ T$
we get (using the identifications $\Lambda ^2V=\mathfrak{so}(V)$)
\begin{eqnarray*}
T(X,Y)&=&\langle T(X,Y),\Theta _{\alpha}\rangle \Theta_{\alpha}
= \langle \mathfrak {T}(X\wedge Y),\Theta_{\alpha}\rangle \Theta_{\alpha}\\
&=& \langle \mathfrak{T}(\Theta _{\alpha}),X\wedge Y)\Theta_{\alpha}
=-\langle T(\Theta _{\alpha})X,Y\rangle\Theta_{\alpha}.
\end{eqnarray*}
The Weizenb\"ock operator is a purely algebraic notion and can be defined for $T$:
$$(\mathcal W^T s)(X_1,...,X_k) = \sum _{i,j} (T(e_j,X_i)s)(X_1,...,e_j,...,X_k)$$
for any tensor $s$ of type $(l,k)$, $k>0$. If $k=0$, we set $\mathcal W ^Ts=0$.
\begin{lemma}For any Riemann-curvature-like tensor $T$ and any  $(0,k)$-tensor $s$ we have
\begin{equation}
\mathcal W ^T s=-\sum_{\alpha} T(\Theta_{\alpha}) (\Theta _{\alpha}s),
\end{equation}
where $T(\Theta_{\alpha})$ acts on $(\Theta_{\alpha}s)$ (and $\Theta_\alpha $ acts on $s$) as   a differentiation.
\end{lemma}

\proof Using the above formula we obtain
\begin{eqnarray*}
(\mathcal W^T s)(X_1,...,X_k) &=& \sum _{i,j} (T(e_j,X_i)s)(X_1,...,e_j,...,X_k)\\
&=& -\sum _{i,j,\alpha} \langle T(\Theta _\alpha)e_j,X_i)(\Theta_\alpha s)(X_1,...,e_j,...,X_k)\\
&=&-\sum_{i,j,\alpha} (\Theta _\alpha s)(X_1,..., \langle T(\Theta_\alpha)e_j,X_i\rangle e_j,...,X_k)\\
&=&\sum _{i,\alpha} (\Theta _\alpha s) (X_1,...,T(\Theta_\alpha)X_i,...,X_k)\\
&=& -\sum _{\alpha}(T(\Theta_\alpha)(\Theta_\alpha s))(X_1,...,X_k).
\end{eqnarray*}
\koniec

\begin{lemma}\label{positive_curvature_operator}
If $T$ is a Riemann-curvature like tensor and  $\mathfrak T$ is the curvature operator for $T$  then
\begin{equation}
\langle \mathcal W ^Ts,s\rangle= \lambda _\alpha |\Theta _\alpha s|^2\end{equation}
 for any $(0,k)$-tensor $s$, where $\Theta _\alpha$ is an orthonormal eigenbasis for $\mathfrak T$ and $\lambda _\alpha$ are corresponding eigenvalues. In particular, if $\mathfrak{T}\ge 0$, then $\langle \mathcal W ^Ts,s\rangle\ge 0$.
\end{lemma}

\proof  Observe first that if $A: V\to V $ is a  skew-symmetric endomorphism, then
$A$  acting on tensors as a differentiation is also skew-symmetric.
 Using the above lemma we now get
\begin{eqnarray*}
\langle \mathcal  W^Ts,s\rangle&=&-\sum _{\alpha} \langle T(\Theta_\alpha)(\Theta_\alpha s),s\rangle\\
&&\ \ \ \ \ \   \sum _\alpha \langle (\Theta _\alpha s),T(\Theta_\alpha)s\rangle= \sum_{\alpha}\lambda _\alpha|\Theta_\alpha s|^2
\end{eqnarray*}
\koniec

 \begin{lemma}\label{about_omega=0}
 {\rm a)} If for a $k$-form $\omega$, where $0<k<n$,  and for every  $A\in \mathfrak{so}(V)$ we have $A\omega=0$, then $\omega=0$.
 \newline
 {\rm b)}
 If for a Riemann-curvature-like tensor $T$  and for every  $A\in \mathfrak{so}(V)$ we have $AT=0$
 then $T$ is a multiple of $R_0$.
  \end{lemma}
  \proof
a) Suppose that $\omega\ne 0$. Let $e_1,...,e_n$ be an orhonormal basis such that $\omega (e_1,...,e_k
)\ne 0$. Take $A\in \mathfrak{so}(V)$ such that $Ae_1=...=Ae_k=0$ and $Ae_{k+1}=e_k$.
Then we get the following contradiction
$$0=(A\omega )(e_1,...,e_{k-1}, e_{k+1}) =-\omega (e_1,..., e_k)\ne 0.$$

b) It is sufficient to observe that if $X,V,W,Z$ are mutually
orthogonal then $T(X,V,W,Z)=0$. First we take three orthogonal
vectors $X,Y,Z$ and $A\in \mathfrak{so}(V)$ such that $AY=0$ and
$AX=Z$. We have $0=(AT)(X,Y,Y,X)=
-T(AX,Y,Y,X)-T(X,Y,Y,AX)=-2T(X,Y,Y,Z)$. Take now $X,Z,W,V$
orthogonal. Then the vectors  $Y=V+W$, $X,Z$ are orthogonal and from
the above formula we get $ T(X,V,W,Z)=-T(X,W,V,Z)$. Finally we
obtain
\begin{eqnarray*}
T(X,V,W,Z)&=&T(W,V,X,Z)-T(W,X,V,Z)= -2T(W,X,V,Z)\\
&&=2T(X,W,V,Z)=-2T(X,V,W,Z)
\end{eqnarray*}
which implies that $T(X,V,W,Z)=0$.\koniec

\subsection{The sectional $\nabla$-curvature}
Let $(g,\nabla)$ be a statistical structure on an $n$-dimensional manifold $M$.
In general, $g(R(X,Y)U, Z)$ is not skew-symmetric for $U,Z$.
Define the following  tensor field of type $(0,4)$
\begin{equation}
  R (U,Z,X,Y)=\frac{1}{2}(g(R(X,Y)Z,U)-g(R(X,Y)U,Z))
\end{equation}
Of course it is skew-symmetric relative to the both pairs of arguments $X,Y$ and
$U,Z$. Since $g(R(X,Y)Z,U)=-g(\overline R(X,Y)U,Z)$ we have
\begin{equation}
 R(U,Z,X,Y)= \frac{1}{2}g(R(X,Y)Z+\overline
R(X,Y)Z,U).\end{equation}
It follows that  the first Bianchi identity holds:
 $$\Xi
_{Z,X,Y} R (U,Z,X,Y)=0,$$ where $\Xi _{Z,X,Y}$ stands for the cyclic
permutation sum relative to $Z,X,Y$. Consequently $R$ is a
Riemann-curvature-like tensor.
We can now define the sectional $\nabla$-curvature of a vector plane
$\pi$ spanned by the orthogonal vectors $e_1, e_2$ by the formula
$$ k(\pi) =k(e_1\wedge e_2)=R(e_1,e_2,e_1,e_2)=
\frac{1}{2}g(R(e_1,e_2)e_2 +\overline R(e_1,e_2)e_2, e_1).$$

We  also have
\begin{equation}
k(X\wedge Y) =\frac{ R (X,Y,X,Y)}{g(X,X)g(Y,Y)-g(X,Y)^2}
\end{equation}
We shall say that  the sectional $\nabla$-curvature is point-wise
constant if for each point $x\in M$ the sectional $\nabla$-curvature
is independent of a plane in $T_xM$ and it is equal to $k(x)$. Of
course, in such a case,  $k(x)$ is a smooth function and, as  in the
case of Levi-Civita connections, one has

$$R(X,Y)Z +\overline R(X,Y)Z =2k(x)\{g(Z,Y)X-g(Z,X)Y\}=:2k(x) R_0(X,Y)Z.$$

Examples of manifolds with constant sectional $\nabla$-curvature are
locally strongly convex equiaffine spheres (in other terminology
relative spheres), where $\nabla$ is the induced connection and $g$
is the affine second fundamental form. If the corresponding  shape
operator is equal to $\lambda \id$, where $\lambda \in \R$, then, by
the Gauss equation,
$$ R(X,Y)Z=\overline R(X,Y)Z=\lambda(g(Y,Z)X-g(X,Z)Y).$$
Hence for  such a sphere the sectional $\nabla$- curvature equals to
$\lambda$.

 In contrast with the classical  sectional curvature, the fact
that the sectional $\nabla$-curvature is point-wise constant does
not imply that it is constant on a connected manifold if the
dimension of the manifold is greater than 2. To see this let us
consider the following example.

\begin{example}
{\rm Take the hypersurface $M$ in $\R ^{n+1}$ given by the equation
$$x_1\cdot ...\cdot x_{n+1}=1$$
for $x_1>0,...,x_{n+1}>0$. It is a locally strongly convex proper
affine sphere (the affine shape operator equals to $\lambda \id$,
where $\lambda\ne 0$) and its Blaschke metric $g$ is flat, i.e. $\hat R=0$.
 Since
$$R(X,Y)Z +\overline R(X,Y)Z= 2\hat R(X,Y)Z +2[K_X,K_Y],$$ we have
$$\lambda =g([K_X,K_Y]Y,X)$$
for any orthonormal pair  of vectors $X,Y$. Take
  $$\tilde\nabla=\hat\nabla +\varphi K. $$
Denote by $\tilde R$ the curvature tensor of $\nabla$ and by
$\overline{\tilde R}$ the curvature tensor of the conjugate
connection $\overline{\tilde\nabla}$. We have
$$\tilde R(X,Y)Z +\overline{\tilde R}(X,Y)Z=2\hat R(X,Y)Z +2\varphi
^2[K_X,K_Y]Z$$ and consequently
$$\frac{1}{2}[g(\tilde R(X,Y)Y,X) +g(\overline{\tilde R}(X,Y)Y,X)]
= \varphi ^2\lambda$$ for orthonormal vectors $X,Y$. If we take a
non-constant function $\varphi$,  we get a structure of point-wise
constant but non-constant $\tilde\nabla$-sectional curvature on $M$.
}\end{example}




Let now $M$ be a locally strongly convex hypersurface with an
equiaffine transversal vector field $\xi$. Let $g$ be the  second
fundamental form,  $S$ -- the corresponding shape operator and
$\nabla$ -- the induced connection. By the Gauss equation we have
for orthonormal $X,Y$
\begin{eqnarray*}
g(R(X,Y)Y+\overline R(X,Y)Y,X)
= g(SX,X) +g(SY,Y).
\end{eqnarray*}
For equiaffine hypersurfaces we can define the
sectional mean curvature. Namely, if $\pi$ is a plane in the tangent
space, then $$k(\pi)= g(Se_1,e_1)+g(Se_2,e_2),$$ where $e_1, e_2$ is
an orthonormal basis of $\pi$. The above considerations show that
the definition is independent of the choice of  an orthonormal
basis. In particular, an equiaffine surface in $\R^3$ of constant
sectional  $\nabla$-curvature is exactly a surface of constant
equiaffine mean curvature. Assume now that the sectional mean
curvature for an equiaffine hypersurface is point-wise constant.
Then we have

\begin{eqnarray*}
g(Y,Z)SX -g(X,Z)SY &+&g(SY,Z)X-g(SX,Z)Y=R(X,Y)Z +\overline R(X,Y)Z\\
&&=2k(g(Y,Z)X -g(X,Z)Y).
\end{eqnarray*}
If $\dim M>2$, then for any $X$ we can take $Z\ne 0$ such that
$g(Z,X)=0$, $g(Z,SX)=0$ and $Y=Z$. We see that $SX$ is a multiple of
$X$. By the second Codazzi equation ($\nabla S$ is symmetric) we
obtain that $S=\lambda \id$, where $\lambda$ is constant if $M$ is
connected. Hence the sectional mean curvature is constant. Roughly
speaking, in the case of equiaffine hypersurfaces Schur's lemma
holds.

Schur's lemma also holds for connections  satisfying the condition
$R=\overline R$.  In the cathegory of Blaschke hypersurfaces the
condition describes affine spheres.

More generally, we have

\begin{lemma}\label{sfery_afiniczne}
Let $M$ be a connected locally strongly convex hypersurface equipped
with a transversal vector field whose induced second fundamental
form is $g$,  the induced connection is $\nabla$ and the induced
shape operator is $S$. The hypersurface is an equiaffine sphere,
that is, $S=\lambda \id$ if and only if $R=\overline R$.
\end{lemma}

\proof Assume that $R=\overline R$. The shape operator is
diagonalizable. Let $\dim M=2$ and $e_1,e_2$ be an orthonormal
basis  of $T_xM$ such that $Se_1=\lambda _1e_1$, $Se_2=\lambda_2
e_2$. By the Gauss equation we have
\begin{equation}\label{from_Gauss}
g(Y,Z)SX-g(X,Z)SY=g(Y,SZ)X-g(X,SZ)Y.\end{equation} Setting $X=e_1$,
$Y=Z=e_2$ we get $\lambda _1=\lambda_2$.

Assume now that $n>2$. Take any $T_xM\ni X\ne 0$ and its orthogonal
complementary space $X^\bot$ in $T_xM$. The mapping
$$X ^\bot \ni W\to g(X,SW)\in \R$$
has kernel of dimension at least $1$. Take $Z=Y\ne 0$ from this
kernel. Using (\ref{from_Gauss}) we obtain that $SX$ is proportional
to $X$, which finishes the proof.\koniec

We have the following second Bianchi identity for the curvature tensor $R+\overline R$
\begin{lemma}
For any statistical structure $(g,\nabla)$ we have
$$\Xi _{U,X,Y}(\hnabla _U(R+\overline R))(X,Y)=
\Xi _{U,X,Y}( K_U(\overline R- R))(X,Y).$$

\end{lemma}
\proof
 We have
 \begin{eqnarray*}&&\Xi _{U,X,Y}(\hnabla _U(R+\overline R))(X,Y))\\
&&\ \ \ \ \ =\Xi _{U,X,Y}((\nabla -K)_UR)(X,Y)+ \Xi _{U,X,Y}((\onabla +K)_U\overline R)(X,Y)\\
&&\ \ \ \ \ \ \ \ \ =\Xi _{U,X,Y}( K_U(\overline R-
R))(X,Y)\end{eqnarray*}\koniec

Thus, if $R=\overline R$ then
$$\Xi _{U,X,Y}(\nabla_UR)(X,Y)= \Xi _{U,X,Y}(\hnabla _UR )(X,Y)=0.$$

Using the second Bianchi identity for $\hnabla R$ one easily gets

\begin{proposition}
Let $M$ be  a connected manifold of dimension greater than 2. If
Let $(g,\nabla)$ be a statistical structure on $M$ with
$R=\overline R$. If the sectional $\nabla$-curvature is point-wise
constant then it is constant.
\end{proposition}

\bigskip

\subsection{ The curvature operator for statistical structures}
If $(g,\nabla)$ is a statistical structure then $T=R+\overline R$ is a Riemann-curvature-like tensor
field and we can apply
the  algebraic  results of section 8.1  to  this tensor field.

Because of the Bianchi identity proved in the last section
exactly in the same way as Theorem 1.2 in \cite{P1} one  can prove
\begin{thm}\label{from_P}
Assume  that for a statistical structure $(g,\nabla)$ we have $R=\overline R$.
The following formula holds
\begin{eqnarray*}
&&(\hnabla^*\hnabla R)(X,Y,Z,W)+\frac{1}{2} (\mathcal W ^{\hat R}R)(X,Y,Z,W)\\
&&\ \ = \frac{1}{2}(\hnabla _X\hnabla ^*R)(Y,Z,W)-\frac{1}{2}(\hnabla _Y\hnabla ^*R)(X,Z,W)\\
&& \ \ \ \  +\frac{1}{2}(\hnabla _Z\hnabla ^*R)(W,X,Y)-\frac{1}{2}(\hnabla_W\hnabla^*R)(Z,X,Y).
\end{eqnarray*}
\end{thm}
We can  now formulate the following version of Tachibana's theorem

  \begin{thm} Let $M$ be a connected compact oriented
manifold and $(g,\nabla)$ be a statistical structure on $M$ such
that $R=\overline R$. If the curvature operator $\hat{\mathfrak
{R}}$ for $\hat R$ is non-negative  and $\div ^{\hnabla}R=0$ then
$\hnabla R=0$. If additionally $\hat{\mathfrak {R}} >0$ at some
point of $M$ then the sectional $\nabla$-curvature is constant.
\end{thm}
\proof It is clear that  $\div ^{\hnabla}R=\nabla ^*R$. By Theorem
\ref{from_P} we now have
$$\hnabla ^*\hnabla R +\frac{1}{2} \mathcal W ^{\hat R}R=0.$$
Consequently
\begin{eqnarray*}
0&=&\int _M (g(\hnabla ^*\hnabla R,R) +\frac{1}{2}g(\mathcal W^{\hat R}R, R))\gnu\\
&&\ \ \ \ \ = \int _M g(\hnabla R,\hnabla R)\gnu +\frac{1}{2}\int _M g(\mathcal W^{\hat R}R, R)\gnu\\
\end{eqnarray*}
By Lemma \ref{positive_curvature_operator}
we obtain
$\hnabla R=0$ and $ \sum _\alpha \lambda _\alpha |\Theta _\alpha R|^2=0$ at each point of $M$.
Therefore, if at some point $x\in M$ the curvature operator $\hat {\mathfrak R}$ is positive, then $\Theta _\alpha R=0$
at this point for all $\alpha$ and consequently for any $A\in \mathfrak{so}(T_xM)$ we have $AR =0$. By Lemma
\ref{about_omega=0} b) we get that $R=\lambda R_0$  at $x$. Since $\hnabla R=0$, the same equality holds at each point of $M$.\koniec

Finally we   observe that  a theorem of  Meyer-Gallot
 holds for trace-free statistical structures

\begin{thm}
Let $M$ be a connected compact oriented  manifold and $(g,\nabla)$
be a trace-free statistical structure on $M$. If the curvature
operator for $R+\overline R$ is non-negative on $M$ then each
harmonic form is parallel relative to $\nabla$, $\onabla$ and
$\hat\nabla$. If moreover the curvature operator is positive at some
point of $M$, then the Betti numbers $b_1(M)=...=b_{n-1}(M)=0$.
\end{thm}

\proof If $\omega$ is a harmonic form, then by
(\ref{Bochner_Weitzenbock_formula_for_s^2}) and Lemma \ref{positive_curvature_operator}
we obtain
\begin{eqnarray*}0&=&\int _M (g(\mathcal W^{R+\overline R}\omega
,\omega) +g(\nabla \omega,\nabla \omega) +g(\onabla
\omega,\onabla\omega ))\gnu \\&=& \int_M \sum _\alpha \lambda
_\alpha |\Theta_\alpha\omega|^2\gnu +\int _M g(\nabla \omega,\nabla
\omega)\gnu +\int_M g(\onabla \omega,\onabla\omega )\gnu.
\end{eqnarray*}
This yields the first assertion. If at some point $x\in M$ all $\lambda
_\alpha
>0$ then we additionally  have  $\Theta _\alpha\omega=0$
at this point.  Since the $\Theta _\alpha$ form a basis for the
space of $\mathfrak{so}(T_xM)$, we have that $A\omega=0$ for all $A\in \mathfrak{so}(T_xM)$.
 Using now Lemma \ref{about_omega=0}  we see that $\omega _x =0$
if the degree of $\omega$ is between $1$ and $n-1$. Since $\omega$
is parallel relative to a connection, it must vanish on the whole of
$M$. \koniec

\end{document}